\documentclass[12pt,reqno,twoside]{amsart}
\usepackage{graphicx}
\usepackage{amsmath}
\usepackage{amssymb}

\usepackage{color}

\newcommand{\ignore}[1]{}

\newcommand{\R}{{\mathbb{R}}}  
\newcommand{\Z}{{\mathbb{Z}}}

\newcommand{\vre}{\varepsilon}

\newcommand\vare{\varepsilon}

\newcommand\CC{\mathbb C}

\newcommand\NN{\mathbb N}
\newcommand\RR{\mathbb R}
\newcommand\ZZ{\mathbb Z}
\newcommand\QQ{\mathbb Q}
\newcommand\TT{\mathbb T}

\newcommand\cH{\mathcal{H}}

\newcommand\vol{\operatorname{vol}}

\newcommand\norm[1]{\left\|#1\right\|}

\newcommand\abs[1]{\left|#1\right|}

\newcommand\set[1]{\left\{{#1}\right\}}

\newtheorem{thm}{Theorem}[section]

\newtheorem{cor}[thm]{Corollary}
\newtheorem{defi}[thm]{Definition}

\numberwithin{equation}{section}

\begin{document}

%
\title[Ergodic theorems and number theory]{Quantitative ergodic theorems and their number-theoretic applications}
%


\date{\today}



\author{Alexander Gorodnik}
\address{School of Mathematics \\ University of Bristol \\ Bristol, U.K.}
\email{a.gorodnik@bristol.ac.uk}

\author{Amos Nevo}
\address{Department of Mathematics, Technion, Israel}
\email{anevo@tx.technion.ac.il}

\begin{abstract} 
We present a survey of ergodic theorems for actions of 
algebraic and arithmetic groups recently established by the authors,  as well as some of their applications. Our approach is based on spectral methods employing the unitary representation theory of the groups involved. This allows the derivation of  ergodic theorems with a rate of
convergence, an  important phenomenon which does not arise in classical ergodic theory.  Quantitative ergodic theorems give rise  to new and previously inaccessible applications, and  we demonstrate the remarkable diversity of  such applications by presenting  several number-theoretic results. 
These include, in particular, general uniform error estimates in lattice points counting problems,
explicit estimates in sifting problems for almost-prime points on symmetric varieties, 
bounds for exponents of intrinsic Diophantine
approximation, and results on fast distribution of dense orbits on homogeneous spaces. 
\end{abstract}

\maketitle
{\small \tableofcontents}

%
%
%

\section{Introduction}
The goals of the present survey are to describe some mean and pointwise ergodic theorems for actions of
algebraic groups and their lattice subgroups, and to give an exposition of  some of the diverse
applications that such ergodic theorems have, particularly to counting, equidistribution and Diophantine
approximation problems in number theory. The results described below  have been recently developed by the authors, and constitute an expansion of the range of applications of spectral methods  in ergodic theory, particularly in the context of dynamics on homogeneous space. The survey is however limited in scope in that it describes only the authors'  specific spectral approach,  and it does not discuss other approaches within  the very large and rapidly expanding field of homogeneous dynamics and its application in number theory.  

In order to situate our approach in the proper context, we start by describing below some elements in the
interesting story of the development  of  the mean ergodic theorem, beginning with its formulation by
von-Neumann in 1932. Let us however note  here very briefly the main new points in our discussion.
First,  the ergodic theorems, which we establish, are developed using spectral methods and representation
theory of non-amenable algebraic groups, subjects left out of most of traditional ergodic theory. Second,
these ergodic theorems provide a rate of convergence to the ergodic limit, a remarkable and most useful
phenomenon that does not arise in classical ergodic theory. 
Third, the rate of convergence allows new and previously inaccessible applications of ergodic theorems to
be developed, using the fact that  the actions of the groups we study are closely connected to a variety
of very natural number-theoretic problems.  

To demonstrate the utility and diversity of such applications, let us briefly describe a selection of
counting, equidistribution,  and Diophantine approximation problems that may at first sight appear
completely unrelated. The main rationale of our exposition in the present survey is the fact, which will
be explained in detail below, that solutions to all of the following problems
depend on applying a suitable ergodic theorem to the case under consideration. 

\subsubsection*{\bf 1. The uniform lattice point counting problem}
Let $\Gamma$ be a discrete subgroup of a locally compact group $G$ such that the space 
$G/\Gamma$ carries a finite $G$-invariant measure.
The lattice point counting problem is to count the
number of elements of $\Gamma$ in an increasing  family of bounded domains $B_t\subset G$, and ideally to
establish an asymptotic formula 
$$
\abs{\Gamma\cap B_t}=\frac{\abs{B_t}}{\abs{G/\Gamma}}+O(\abs{B_t}^\kappa)
$$
with $\kappa < 1$ as small as possible. 
Some questions that arise here are, for example, 
what is the asymptotics of the number of unimodular integral matrices $A$
lying in the norm ball $\|A\|<t$ as $t\to\infty$, 
or in more general domains in the matrix space ? 
  What is the asymptotics for the number of such matrices $A$ satisfying a congruence condition
  $A=A_0\mod q$?
More generally, what is the asymptotics of the number of unimodular matrices 
with entries in the ring of
integers of an algebraic number field, which are of bounded height? 


Counting lattice points is in fact the most basic of all mean ergodic theorems that we consider. This point of view that will be explained in full in \S 5 below, where 
we will describe a  general solution to the lattice point counting problem. 

\subsubsection*{\bf 2.  Almost prime points on algebraic varieties}

As we shall see  below, establishing robustness and uniformity properties of the solution to the lattice
point counting problem is crucial in several important applications, which depend on answering  more
refined counting questions.  Of those, let us mention the problem of
existence of {\it almost prime} integral points on homogeneous varieties, raised in \cite{ns}.  Consider for example
the set of symmetric integral positive-definite matrices of fixed determinant, and call such a matrix prime if all
its entries are prime numbers.  How large is the set of prime points in the set of integral points?
Less ambitiously, how large is the set of  integral points all of whose entries are $r$-almost primes 
(i.e., products of at most $r$ prime factors) ?

This question has been the subject of intensive activity in recent years, and we refer  to \cite{bgs2}
and the recent surveys \cite{Lu} and \cite{kont} for further discussion. In \S 6 below we will utilise
uniformity in the solution of the lattice point counting problem to 
establish an asymptotic lower bound of the right magnitude
for the number of $r$-prime points lying on symmetric varieties.



\subsubsection*{\bf 3.  Diophantine approximation on algebraic varieties.} 
Consider the rational ellipsoid  $\set{(x,y,z)\,;\, ax^2+by^2+cz^2=d}$, $a,b,c,d\in \NN_+$.
Assuming that the rational points on the ellipsoid are dense, 
is it possible to establish a rate of approximation of a general points on the ellipsoid by rational
points ? Is it possible to establish a rate  of approximation by rational points satisfying  additional
integrality constraints ? 

 These problems, in the much more general context of homogeneous algebraic varieties,  were raised in
 S. Lang's 1965 Report on  Diophantine Approximation \cite{L65}. We will state this problem precisely in
 \S 7 and describe a general solution, obtained jointly with A. Ghosh,   in the case where the variety
 is homogeneous under an action a simple algebraic group. For the rational ellipsoids we establish a rate of approximation by rational points satisfying integrality constrains  that is best possible in a natural sense. 

\subsubsection*{\bf 4. Equidistribution of typical orbits  on homogeneous spaces.} 

Consider the one-sheeted hyperboloid $\set{(x,y,z)\,;\, ax^2+by^2-cz^2=d}$, $a,b,c,d\in \NN_+$,
which is also known as de-Sitter space. The group of integral matrices preserving the corresponding
quadratic form naturally acts on this space, and typical orbits are dense, but does a typical orbit
have a limiting distribution in a suitable sense ? 

 This problem,  for general lattices in the isometry group of a Lorentz form acting on de-Sitter space of arbitrary dimension, was  raised by V. Arnold in 1985 \cite{a}. A precise  quantitative solution for general homogeneous varieties will be described in \S 9.   

\subsubsection*{\bf 5.  Equidistribution in algebraic number fields.} The group  $\Gamma= {\rm
  SL}_2(\ZZ[\sqrt{d}])$, where $d >0$ is not a square,  is a dense subgroup of $G={\rm SL}_2(\RR)$.
Consider elements of $\Gamma$ of bounded norm : $\Gamma_t=\{\gamma\in\Gamma:\, \log(\|\gamma\|^2+\|\bar
\gamma\|^2)\le t\}$, where $\bar\gamma$ denotes the Galois conjugate. 
Do these sets become equidistributed in  $G$ in a suitable sense ? 

We note that the action of $\Gamma$ on $G$ by translations is isometric, and the natural problem is to prove equidistribution for {\it every single orbit}, i.e. every coset of $\Gamma$. One possible such statement is a ratio equidistribution theorem, establishing that for nice domains $\Omega_1,\Omega_2\subset G$, 
$$\frac{\abs{ \Gamma_t \cap \Omega_1}}{\abs{  \Gamma_t \cap  \Omega_2}}\to \frac{|\Omega_1|}{|\Omega_2|}\quad\hbox{as $t\to\infty$} \,,$$
with the convergence taking place at a definite rate. Such a result will be stated precisely and explained in \S 9. 

We remark that the problem of  establishing a ratio equidistribution theorem for dense subgroups acting
by isometries on the group or its homogeneous spaces was raised originally by Kazhdan \cite{k}, in the
case of a dense group of Euclidean isometries acting on the Euclidean plane.  We refer to \cite{v}
for a recent progress on this problem.

\subsubsection*{\bf 6. Ergodic theorems for lattice subgroups}  

Finally, let us mention some natural questions which constitute direct generalizations of the fundamental convergence problems of classical ergodic theory.
 Taking the  group $\Gamma=\hbox{SL}_n(\ZZ)$ as an example, consider its linear action 
on the torus $\TT^n=\mathbb{R}^n/\mathbb{Z}^n$. Does the mean ergodic theorem hold for the action of $\Gamma$ on $\TT^n$  ? Specifically, let 
$\Gamma_t =\set{\gamma\in \Gamma\,;\, \log \norm{\gamma}\le t}$, and for $f\in L^p(\TT^n)$ and $x\in\TT^n$,
consider the averages
$$
\frac{1}{|\Gamma_t|}\sum_{\gamma\in\Gamma_t}f(\gamma^{-1}x).
$$
Do these averages converge to the space average $\int_{\TT^n} f(x)\, dx$ in $L^p$? 
Does convergence hold almost everywhere ?
If so, is there a rate of convergence to the space average, in $L^p$ and almost everywhere ?

 
 
These questions also arise for other interesting actions of $\hbox{SL}_n(\ZZ)$ 
(for instance, for the action on the space of unimodular lattices or for the action on its profinite
completion), and they
are all resolved by a general pointwise and mean ergodic theorem for the averages on norm balls in $\hbox{SL}_n(\Z)$, a result that will be formulated more generally and precisely in \S 8 and explained there. 
\vskip0.3truein

To conclude the introduction, let us emphasize that it is ultimately the existence of an explicit rate of
convergence in the ergodic theorems that is responsible for solution
of the above problems.
  We will therefore give in our discussion high priority to  the problem of establishing  an explicit  rate of convergence of the  averages  we consider to the space average in measure-preserving dynamical systems.

To situate our approach in its proper context, let us now turn to describe some elements in the development of the mean ergodic theorem.  We will also mention some of the milestones in the applications of spectral methods in ergodic theory. As we shall see, an interesting  feature of this story is that the success of the classical methods of ergodic theory is inextricably tied to both the impossibility of establishing a rate of convergence in the ergodic theorems, and the impossibility of extending the classical methods to classes of groups such as simple algebraic groups or their lattice subgroups. 

\section{Averaging in  ergodic theory }

\subsection{The classical mean ergodic theorem}

 The basic objects of study in classical Hamiltonian mechanics consist of  a compact Riemannian manifold
 $M$, which is called the phase space, together with  a divergence-free vector field on $M$.  We denote by $\mu$ the Riemannian volume on $M$, normalised to be a probability measure. 
 The integral curves of  the vector field give rise to a one-parameter group of volume-preserving
 transformations  $a_s : M\to M$, $s\in \RR$, which describes the time evolution in phase space.
For a function $f : M \to \RR$, the  time averages along orbits are defined by 
$$\mathcal{A}_tf(x)=\frac1t \int_0^t f(a_s x)\, ds.$$
One of  the most fundamental questions in classical dynamics  is the elucidation of their long term
behavior. Two important questions that arise are : 
\begin{itemize}
\item do the time averages $\mathcal{A}_tf$ converge at all, and if so, what is their limit?
\item how does the limit depend on the initial state $x$ and the function $f$?
\end{itemize}

%
%
The solution proposed by von Neumann  to this problem was to consider for every measure-preserving map
$a_s : M \to M$ an associated unitary transformation, the Koopman operator $U_s : L^2(M)\to L^2(M)$,
given by $U_s f (x)=f(a_s x)$. Thus  the measure-preserving flow gives rise to a one-parameter group of
unitary transformations. This made it possible for von-Neumann to apply his recently developed spectral
theory to prove the mean ergodic theorem \cite{vN31}: 
\begin{thm}[von Neumann]\label{th:vn}
Let $(X,\mu)$ be a (standard) probability space  with a one-parameter flow of
measure-preserving transformations  $a_s : X\to X$,  $s\in \RR$.  
Then for any $f\in L^2(X)$, 
$$
\mathcal{A}_tf\to \mathcal{P}f\quad \hbox{as $t\to \infty$},
$$
in $L^2$-norm, where $\mathcal{P}$ denotes the 
orthogonal projection on the space of functions invariant under the flow.
\end{thm}
 In particular, consider the case where the  space of invariant functions in $L^2(M)$ consists of constant functions
 only, or equivalently, every invariant measurable set has measure zero or one. Flows satisfying this
 property are called \emph{ergodic}. Then the orthogonal projection of $f$ to the space of constant
 functions is given by the ``space average'' $\int_X f\, d\mu$. Thus, for ergodic flows the time
 averages converge in $L^2$-norm to the space average:
$$
\mathcal{A}_t f \to \int_X f\, d\mu\quad\hbox{as $t\to \infty$},
$$
verifying the ``ergodic hypothesis''. 

Several years later, F. Riesz's \cite{Rie} gave an elegant elementary  proof of the mean ergodic theorem,
which  became very influential in the subsequent development of ergodic theory. We will indicate it here in full, in order to 
highlight its main idea. 
Let $U_s :\cH\to \cH$, $s\in \mathbb{R}$, be a one-parameter group of unitary operators on a Hilbert space $\cH$
and $\mathcal{P}$ denote the projection on the space of vectors invariant under the group. 
We will show that the averaging operators $\mathcal{A}_t=\frac1t\int_0^t U_s \,ds$ 
satisfy $\mathcal{A}_t f \to \mathcal{P} f$ for every $f\in\mathcal{H}$
by breaking $\cH$ up  to two orthogonal complementary subspaces and arguing separately in each of them. 
First, we apply $\mathcal{A}_t$  to vectors $f$ of the form $f=U_{t_0} h-h$. Then
\begin{align*}
\mathcal{A}_t f &= \mathcal{A}_t(U_{t_0} h-h)= \frac1t\left(\int_{t_0}^{t+t_0} U_s h\,ds -\int_0^t U_s h \,ds\right)\\
&=\frac1t \left( \int_t^{t+t_0}-\int_0^{t_0}\right)U_s h \,ds,\nonumber
\end{align*} 
and therefore 
\begin{equation*}
\norm{\mathcal{A}_t f}\le \frac{2t_0}{t}\norm{h}\rightarrow 0\quad\hbox{as $t\to \infty$.}
\end{equation*}
Since $f$ is orthogonal to every invariant vector, namely $\mathcal{P} f=0$, the mean ergodic theorem holds for $f$ as stated, hence also for the closed linear span of such $f$'s.
Obviously, the mean ergodic theorem holds in the space of invariant vectors. The proof is therefore
complete upon observing that the  span of $U_{t_0}h-h$ with $t_0\in \RR$ and $h\in \cH$ is dense in the orthogonal complement of the space of invariant vectors.

Notice that the main ingredient of the above proof is the property of asymptotic invariance (under
translation) of the intervals $[0,t]$, namely the fact that the measure of the symmetric difference of
$[0,t]$ and $[0,t]+t_0$ divided by the measure of $[0,t]$ converges to zero at $t\to \infty$.
This crucial property of intervals on the real line admits a straightforward generalisation. 
Let $G$ be a locally compact second countable group $G$ equipped with a left-invariant Haar measure.
Define a one-parameter family of sets $F_t\subset G$ of positive finite measure to be
\emph{asymptotically invariant} if it  satisfies 
\begin{equation}\label{eq:folner}
\frac{\abs{F_t g\triangle F_t}}{\abs{F_t}}\rightarrow 0\quad\text{ as }t\to \infty.
\end{equation}
uniformly for $g$ in compact sets of $G$.
A group admitting such an asymptotically invariant family (also called a \emph{F\o lner family}) is
called \emph{amenable}. 

Given a measure-preserving action of the group $G$ on a probability space $(X,\mu)$,
 define for $g\in G$ the unitary Koopman operators
$$
\pi_X(g):L^2(X)\to L^2(X): f\mapsto f(g^{-1}x).
$$
It is straightforward to adapt  Riesz's argument and  conclude that amenable
groups satisfy  the mean ergodic theorem.
Namely, denoting by $\beta_t$ the uniform probability measures supported on 
an asymptotically invariant (F\o lner) family $F_t$ and by $\pi_X(\beta_t)$ 
the corresponding averaging operators in $L^2(X)$, we conclude that 
$$
\pi_X(\beta_t) f
=\frac{1}{\abs{F_t}} \int_{F_t} \pi_X(g)f\, dg\to \mathcal{P}f\quad\hbox{as $t\to \infty$,}
$$
in $L^2$-norm.

\subsection{Equidistribution and pointwise convergence}
So far we have discussed the convergence of the averages $\pi_X(\beta_t)f$ in $L^2$-norm. 
The question as to the convergence $\pi_X(\beta_t) f(x)$ for individual points $x\in X$ is of fundamental importance as well. 

We  begin by describing the best-case scenario, namely equidistribution, and its characterisation in amenable dynamics.
Let $X$ be a compact metric space on which $G$ acts continuously, and 
let  $F_t \subset G$ be a family of F\o lner sets that define the corresponding 
averaging operators $\pi_X(\beta_t)$.
An extreme possibility is that for every continuous function $f$ and for { every} point $x\in X$, the ``time averages'' 
$\pi_X(\beta_t) f(x)$ converge uniformly 
to a limit $m(f)$, which is independent of $x$. In that case, the functional $f\mapsto m(f)$ defines a 
a probability measure $m$ on $X$, which is invariant under $G$. Furthermore $m$ is then the unique $G$-invariant probability measure on $X$, and the system is called \emph{uniquely ergodic}. Conversely, when there exists a unique  $G$-invariant probability measure on $X$, the ``time averages'' of $f$ converge uniformly to the space average $\int_X f\, dm$. 

In general the situation is far more complicated, and the limit of the averages $\pi_X(\beta_t)f(x)$ 
depends very sensitively on the initial point $x$. Typically there exist two 
(or even uncountably many) mutually singular  
invariant probability measures $\mu_1$ and $\mu_2$ on $X$, such that the following holds. 
The time averages $\pi_X(\beta_t) f(x)$ of a continuous function $f$ 
converge on a set $X_1$ of full $\mu_1$-measure to $\int_X f d\mu_1$, and on a set $X_2$ of full
$\mu_2$-measure to $\int_X f d\mu_2$, with the sets $X_1$ and $X_2$ disjoint, and $\int_X f d\mu_1\neq \int_X f
d\mu_2$. 
The question of pointwise convergence for one-parameter flows is addressed by
Birkhoff's pointwise ergodic theorem \cite{Bi}:

\begin{thm}[Birkhoff]\label{th:bf}
Let $(X,\mu)$ be a (standard) probability space  with a one-parameter flow of
measure-preserving transformations  $a_s : X\to X$,  $s\in \RR$.  
Then for any $f\in L^1(X)$, the time averages $\frac{1}{t}\int_0^t f(a_s x)\, ds$
converge for $\mu$-almost every $x\in X$ to $\mathcal{P}f(x)$.
\end{thm}


An elementary proof of Birkhoff's theorem which relies on asymptotic invariance of intervals as well as their order structure was given by \cite{KW}.  So far, this proof has not been generalized to other groups, including $\ZZ^2$. The most successful approach to the proof of poinwise ergodic theorems relies on controlling the maximum deviation between the time averages and the space average throughout  the entire sampling process. Thus we consider the following fundamental quantity, called  the \emph{maximal function} associated to the flow :
$$
\mathcal{M}f=\sup_{t> 0} \abs{\pi_X(\beta_t) f-\mathcal{P}f}.
$$ 
In order to establish the pointwise convergence, 
we seek to prove the following two properties of the flow:
\begin{itemize}
\item[I.] Existence of a dense set of functions  $f\in L^1(X)$ for which $\pi_X(\beta_t) f(x)$ converges for almost
  every $x\in X$. 
\item[II.] The maximum error in the measurement performed by the time averages is controlled  by the
  \emph{strong maximal inequality} in $L^2(X)$:
$$
\norm{\mathcal{M}f}_2\le \hbox{const}\cdot \norm{f}_2,
$$ 
or  by the \emph{weak maximal inequality} in $L^1(X)$:
$$
\mu(\set{ \mathcal{M}f > r})\le \hbox{const} \cdot r^{-1} \norm{f}_1 .
$$
\end{itemize}
Once these properties are established, it follows by a classical approximation argument known as the Banach principle  that
the time averages converge almost everywhere for every $f\in L^1(X)$.
 
The proof of the first property proceeds by using 
asymptotic invariance of a F\o lner family, which implies that
pointwise convergence holds for the family of functions $\pi_X(g)f-f$
with $f\in L^\infty(X)$ and $g\in G$. This is 
a straightforward adaptation of Riesz's argument described above. 

The proof of the maximal inequality utilizes two separate ideas :  the \emph{transference principle} and \emph{covering arguments}. The transference principle asserts 
that the maximal inequality for any family of operators $\pi_X(\beta_t)$ acting on a general space $(X,\mu)$ 
can be reduced to the maximal inequality for the convolution operators defined by $\beta_t$  on the group $G$.
This principle was initiated by Wiener \cite{Wi}, and generalized by Calderon \cite{Ca}
and Coifman--Weiss \cite{CW}. Finally, the investigation of 
the convolution operators exploits \emph{geometric covering properties}
of the translates of sets in a F\o lner family $F_t$.
This approach originated with Wiener \cite{Wi}, who has generalised  the Hardy-Littlewood  maximal
inequality from the real line to Euclidean space by introducing an elegant and influential covering
argument based on the \emph{doubling property}  
$$
|B_{2t}|\le \hbox{const}\cdot |B_t|
$$
for volumes of the Euclidean balls $B_t$. 

The pointwise ergodic theorem was subsequently established, in turn, for groups with
polynomial volume growth, exponential solvable Lie groups, and connected amenable Lie groups. This was
carried out by Templeman \cite{Te1}, Emerson--Greenleaf \cite{EG}, Coifman--Weiss \cite{CW} and others.   
It should be pointed out, however,  that a F\o lner sequence may fail to satisfy the pointwise ergodic
theorem, even for the group $\ZZ$. The pointwise ergodic theorem for general amenable groups requires a
regularity assumption stronger than asymptotic invariance,  namely temperedness, and was established by
Lindenstrauss in \cite{Li}. We refer to  \cite{N3} for a recent survey and to \cite{AAB} for a detailed account of the ergodic theory of amenable groups.


\section{From amenable to non-amenable groups}

The proofs of the ergodic theorems indicated in the previous section were founded upon 
 the almost invariance property (\ref{eq:folner}) of the F\o lner families. 
How then  can one proceed  when the group is not amenable, so that there are no asymptotically invariant families of 
sets at all?  

Before describing our approach, let us note one more consequence of the existence of asymptotically
invariant sequences of sets in the group.
For any properly ergodic (i.e., not transitive) action of a countable amenable group $G$, 
one can show that the space $X$ admits a non-trivial 
asymptotically invariant sequence of sets,  namely, a sequences of measurable sets $A_n \subset X$, such that 
$$
0<\liminf_{n\to \infty} \mu(A_n)\le \limsup_{n\to \infty} \mu(A_n) <1\quad
\hbox{and}\quad \mu(gA_n\triangle A_n)\to 0
$$
uniformly for $g$ in compact subsets of $G$.
It then follows
easily that there exists an \emph{asymptotically invariant sequence} of functions $f_n\in L^2(X)$ 
with zero integral such that 
\begin{equation}\label{eq:ai}
\norm{f_n}=1\quad\hbox{and}\quad \norm{\pi_X(g)f_n-f_n}\to 0
\end{equation}
uniformly for $g$ in compact subsets of $G$.
In fact, the existence of an asymptotically invariant sequence in every property ergodic action
characterises amenable groups (see \cite{Sch1,JR}). 

We can therefore conclude 
 that a countable group $G$ is amenable if and only if in every properly ergodic action,
the averaging operators satisfy 
\begin{equation}\label{eq:n1}
\norm{\pi_X(\beta)|_{L^2_0(X)}}=1
\end{equation}
for every probability measure $\beta$ on $G$, where
$L^2_0(X)$ denotes the orthogonal complement to the space of constant functions.
In particular, it follows that for a properly ergodic action of a countable amenable group, \emph{no uniform rate of convergence} 
in the mean ergodic theorem can be established, namely the norm of the averaging operators does not decay at all.  

\subsection{Spectral gap property}\label{sec:ug}

 Our discussion so far has lead us to the following  observation :  when $G$ is a countable non-amenable group, at least some actions $X$ have the property that the operators $\pi_X(\beta)$ are strict contractions on $L^2_0(X)$, namely $\norm{\pi_X(\beta)} _{L^2_0}< 1$, for all  absolutely continuous generating probability measure $\beta$ on $G$. 
  We can therefore conclude that a non-amenable countable group $G$ necessarily has actions which admit a spectral gap,  in the sense of the following general definition.

\begin{defi}  
{\rm 
A measure-preserving ergodic action of a group $G$ 
on a probability space $(X,\mu)$ has a \emph{spectral gap}  if one of the following two equivalent conditions holds:
\begin{itemize}
\item 
There is no asymptotically invariant sequences of functions  as in (\ref{eq:ai}).
\item  
There exists an absolutely continuous symmetric probability measure $\beta$ whose support generates $G$ such that  
\begin{equation}\label{eq:sg}
\norm{\pi_X(\beta)|_{L^2_0(X)}} < 1.
\end{equation} 
\end{itemize}
}
\end{defi}

We remark that it follows that for actions with spectral gap the estimate 
(\ref{eq:sg}) holds 
for every absolutely continuous probability measure $\beta$ on $G$ whose support generates $G$
as a semigroup.

Kazhdan  \cite{k} has made the truly fundamental discovery that  many groups satisfy an extremely general
form of the spectral gap condition. Not only it is satisfied for all ergodic probability preserving actions of the group,  
an analogous property holds in fact 
for all unitary representations without invariant unit vectors.

\begin{defi}
{\rm 
A group $G$ has \emph{Kazhdan property (T)} if one of the following equivalent conditions holds:
\begin{itemize}
\item Every unitary representation $\pi$ of $G$ which has no 
invariant unit  vectors also has no sequences of asymptotically invariant unit vectors as in (\ref{eq:ai}).
\item For some absolutely continuous symmetric probability measure $\beta$ on $G$ whose support generates $G$, 
\begin{equation}\label{eq:prop_T}
\sup_\pi \|\pi(\beta)\|< 1,
\end{equation}
where the supremum is taken over all unitary representations of $G$ that have 
no invariant unit vectors.
\end{itemize}
}
\end{defi}

Examples of groups satisfying property (T) abound. For instance, 
the groups  
 $\hbox{SL}_n(\RR)$, $\hbox{SL}_n(\CC)$ and $\hbox{Sp}_n(\RR)$ with $n\ge 3$
 have this property.
Among simple  Lie groups, only  the isometry groups of real hyperbolic spaces and the isometry groups of
complex hyperbolic spaces do not have property (T) (the case of rank-one groups is dealt by \cite{Kos}). Similarly, simple algebraic groups over
non-Archemedian local fields all have property (T), unless their split rank is one. We refer to the monograph
\cite{kazhdan_book} for an extensive discussion of this property. 

The isometry groups of real and complex hyperbolic spaces (and algebraic group of split rank $1$) do in fact have probability preserving actions admitting an asymptotically invariant sequence of unit vectors, namely actions which  do not have a spectral gap. Nevertheless, many of the most important and natural actions of these groups do have a spectral gap, although it may vary from action to action, and it does not obey the remarkable uniform estimate provided by property $T$. The problem of establishing uniform spectral estimates for certain families of actions (or families of unitary representations) of these groups is a problem of foundational importance, which has been studied extensively. Let us now turn to describe some of the results and the challenges of establishing such spectral estimates uniformly, and then explain the utility of such estimates in proving ergodic theorems with rates of convergence for actions of semisimple algebraic groups.

\subsection{Selberg's property and property  $\tau$}
The phenomenon where a family of actions with increasing complexity 
satisfies the spectral gap property uniformly, namely the estimate in (\ref{eq:sg}) is uniform,  is termed \emph{property $(\tau)$}.

The first instance of the property $(\tau)$ was discovered by A. Selberg in \cite{selberg}
for congruence subgroups of $\hbox{SL}_2(\Z)\subset \hbox{SL}_2(\mathbb{R})$.
Let us consider a sequence of covers of finite area hyperbolic surfaces
$$
\mathbb{H}^2\rightarrow \cdots\rightarrow S_n\rightarrow \cdots \rightarrow S_2\rightarrow 
S_1=\hbox{SL}_2(\mathbb{Z})\backslash \mathbb{H}^2.
$$
We denote by $\lambda(S_n)$ the bottom of the spectrum 
(excluding the zero eigenvalue) of the Laplace operator
$\Delta=y^2\left(\frac{\partial^2}{\partial x^2}+ \frac{\partial^2}{\partial y^2}\right)$
on $S_n$. One might expect that when the covers $S_n$ grow in a regular fashion,
they ``approximate'' the hyperbolic space $\mathbb{H}^2$, and $\lambda(S_n)$
approaches $\lambda(\mathbb{H}^2)=\frac{1}{4}$.
In particular, $\lambda(S_n)$ is expected to stay bounded away from zero.
Here the assumption that covers grow in a regular fashion is essential, and 
A. Selberg constructed an example of a sequence of covers $S_n$ such that 
$\lambda(S_n)\to 0$ as $n\to\infty$. On the other hand,
for the surfaces
$S_\ell=\Gamma(\ell)\backslash \mathbb{H}^2$, where $\Gamma(\ell)=\{\gamma\in \hbox{SL}_2(\mathbb{Z}):\, \gamma=\hbox{I }\,\hbox{mod}\, \ell\}$ denotes the congruence subgroup 
of $\hbox{SL}_2(\mathbb{Z})$ of level $\ell$, he showed that
\begin{equation}\label{eq:selberg_1}
\lambda(S_\ell)\ge \frac{3}{16}.
\end{equation}
Moreover, A. Selberg conjectured that for the congruence covers,
\begin{equation}\label{eq:selberg_2}
\lambda(S_\ell)\ge \lambda(\mathbb{H}^2)=\frac{1}{4}.
\end{equation}
Although estimate (\ref{eq:selberg_1}) has been improved in 
\cite{gj,lps,i_s,ks}, the Selberg conjecture (\ref{eq:selberg_2}) is still  open.
We refer to the surveys \cite{s0,s1} for a more detailed discussion. 
Nowadays this conjecture is understood as a special case of the generalised
Ramanujan conjectures, discussed in \cite{s,bb}, which are most conveniently stated in the representation-theoretic terms. 
These conjectures also include the original Ramanujan conjecture
concerning the bounds on Fourier coefficients of holomorphic modular forms,
which was solved by P. Deligne.

More generally, let ${\rm G}\subset \hbox{GL}_d(\mathbb{C})$ be a semisimple algebraic group defined over $\mathbb{Q}$,
$\Gamma={\rm G}(\mathbb{Z})$, and let 
$$
\Gamma(\ell)=\{\gamma\in \Gamma:\, \gamma=\hbox{I }\,\hbox{mod}\, \ell\}
$$
denote the congruence subgroups. We denote by 
$$
\mathcal{H}_{\rm G}(\ell)=L_0^2({\rm G}(\mathbb{R})/\Gamma(\ell))
$$ the space of square-integrable functions
with zero  integral on the space ${\rm G}(\mathbb{R})/\Gamma(\ell)$.
We consider the family of unitary representation
of ${\rm G}(\mathbb{R})$ on $\mathcal{H}_{\rm G}(\ell)$ defined by
\begin{equation}\label{eq:pi_q}
\pi_\ell(g)\phi(x)=\phi(g^{-1}x),\quad \phi\in \mathcal{H}_{\rm G}(\ell).
\end{equation}
The generalised Ramanujan conjectures seeks to describe the decomposition
of the representations $\pi_\ell$ (and more general automorphic representations) into
irreducible constituents. We refer to \cite{s,bb} for a comprehensive survey of this subject,
but what would be of crucial importance  for our purposes is the uniform spectral gap property
generalizing (\ref{eq:selberg_1}). In the representation-theoretic language, this property 
can be reformulated in terms of isolation of the trivial representation of ${\rm G}(\mathbb{R})$
from the other irreducible representations appearing in $\pi_\ell$ for all $\ell\ge 1$
or in terms of estimates on the integrability exponents of 
the representations $\pi_\ell$, which we proceed to define. 

If $\Pi$ is a family of unitary representations,
we say that $\Pi$ has {\it property $(\tau)$} if (\ref{eq:prop_T}) holds with
$\pi\in\Pi$. It is of crucial for many
number-theoretic applications to know the uniform spectral gap property
for the family of automorphic representations (\ref{eq:pi_q}).
This problem has a long history and is closely connected with works on Langlands
functoriality conjectures. In particular, for other forms of the group $\hbox{SL}_2$,
property $(\tau)$ follows from Gelbart--Jacquet correspondence \cite{gj}.
An important milestone was reached by M.~Burger and P.~Sarnak \cite{bs} who showed that 
it is sufficient to verify property $(\tau)$ for certain proper subgroups of 
the ambient group. Combined with previous work, this established 
property $(\tau)$ for all 
semisimple $\mathbb{Q}$-simple simply connected algebraic groups
except the unitary groups. Finally, the case of unitary groups
was settled by L.~Clozel in \cite{clo}.

In the case of semisimple Lie groups, it is convenient to measure the quality of spectral gap 
in terms of integrability exponents
that were introduced by Cowling \cite{co}, and Howe and Moore \cite{hm}, as follows.
Let $\pi$ be a unitary representation of a locally compact group $G$
on a Hilbert space $\mathcal{H}$. For $v,w\in \mathcal{H}$,
we denote by $c_{v,w}(g)=\left<\pi(g)v,w\right>$ the corresponding matrix coefficient,
which is a  continuous bounded function on $G$. We define the {\it integrability exponent} of the representation $\pi$ by
\begin{equation}\label{eq:integ}
q(\pi)=\inf\left\{ q>0: 
\, c_{v,w}\in L^q(G)
\hbox{ for all $v,w$ in a dense subset of $\mathcal{H}$}
\right\}.
\end{equation}
It was shown in \cite{bw,co} that when $G$ is a connected simple Lie group with property $T$,
\begin{equation}\label{eq:prop_T2}
\sup_\pi q(\pi)<\infty,
\end{equation}
where the supremum is taken over all unitary representations of $G$ that have 
no nonzero invariant vectors. The value of the supremum has been 
estimated in \cite{ht,li,lz,o,n}. 
The estimate (\ref{eq:prop_T2}) can be considered as a quantitative version of property (T).
We also note the important fact that the works on property $(\tau)$ for automorphic representations also
lead to explicit estimates on $\sup_{\ell\ge 1} q(\pi_\ell)$.

\section{Spectral gaps and ergodic theorems}

\subsection{Quantitative mean and pointwise theorems for $G$-actions}
We now turn to discuss the applications of the spectral estimates described in the previous section in the proofs of ergodic theorems. 
Let  $(X,\mu)$ be a (standard) probability space with a measure-preserving action of 
a locally compact second countable $G$ equipped with left-invariant Haar measure.
For a family of measurable sets $B_t$ of $G$ with finite positive measure, we denote by $\beta_t$
the uniform probability measures supported on $B_t$. We aim to prove ergodic theorems for the corresponding averaging operators $\pi_X(\beta_t)$ defined by
\begin{equation}\label{eq:average0}
\pi_X(\beta_t)f(x)=\frac{1}{|B_t|}\int_{B_t} f(g^{-1}x)\,dg,\quad f\in L^p(X).
\end{equation}



In the rest of this section,
 we discuss quantitative ergodic theorems for connected simple Lie groups $G$
(for instance, $G=\hbox{SL}_n(\mathbb{R})$).
%
In order to prove 
the quantitative mean ergodic theorem for the averages $\pi_X(\beta_t)$, we need to estimate
the norms of the operators $\pi_X(\beta_t)|_{L^2_0(X)}$, and this is achieved  using
the method of ``spectral transfer'',  originating in \cite{N2}.
The methods consists of two separate steps. First, it is possible to reduce the norm estimate of $\pi(\beta_t)$ in a general  representation $\pi$, 
to a norm estimate for the convolution operator $\lambda_G(\beta_t)$ for the regular representation
$\lambda_G$ of 
$G$ acting on $L^2(G)$. This serves as an analogue in this context, of the transference principle for amenable groups.  Second, it is then necessary to estimate the norm of the convolution operators on $G$. 
\ignore{
When the action has a spectral gap, the mean ergodic theorem can be greatly improved by establishing a rate of convergence, and we now explain a general method of "spectral transfer"  (originating in \cite{N2}) for doing so. 
Our goal is to  give a general estimate of $\norm{\pi_X(\beta_t)}$ as an operator on $L^2_0(X)$, 
and the idea of "spectral transfer" is to reduces the norm estimate 
of the operator $\pi(\beta)$ in a {\it general  representation} $\pi$, to a norm estimate for the convolution operator $\lambda_G(\beta)$ in the {\it regular  representation} of $G$.  Thus it serves as an analog in this context, of the transference principle for amenable groups.  
}

Put in quantitative form, the principle of spectral transfer is based on the following two fundamental facts:
 \begin{itemize}
 \item 
If the action of $G$ on $(X,\mu)$ has spectral gap, then the representation $\pi_X|_{L^2_0(X)}$
has finite integrability exponent, and it follows that 
sufficiently high tensor power representation $(\pi_X|_{L^2_0(X)})^{\otimes n}$
embeds as a subrepresentation of the regular representation $\lambda_G$ of $G$.
It then follows from Jensen's inequality that for even $n$,
$$
\norm{\pi_X(\beta_t)|_{L^2_0(X)}}\le \norm{\lambda_G(\beta_t)}^{1/n}.
$$

\item  The convolution operator $\lambda_G(\beta_t)$ 
in the regular representation obey a remarkable convolution inequality, known as \emph{the Kunze--Stein
phenomenon}:  for $\beta\in L^r(G)$ with $1 \le r < 2$ and $f\in L^2(G)$, 
$$
\norm{\beta\ast f}_{2}\le c_r \norm{\beta }_{r}\norm{f}_2$$
with uniform $c_r>0$. 
This convolution inequality for $\hbox{SL}_2(\RR)$ is due to Kunze and Stein \cite{ks0} and in general to
Cowling \cite{co0}. The Kunze--Stein phenomenon therefore implies that
$$
\norm{\lambda_G(\beta_t)}\le c_r \left\|\frac{\chi_{B_t}}{|B_t|}\right\|_r =c_r |B_t|^{-1+1/r}.$$
\end{itemize}
This leads to  an estimate for $\left\|\pi_X(\beta_t)|_{L^2_0(X)}\right\|$, and
using an interpolation argument, one extends this estimate for the action of $G$ on $L^p(X)$, $ 1 < p < \infty$.
Hence, we deduce the following quantitative mean ergodic theorem, originated in \cite{N2}:

\begin{thm}\label{spectral transfer}  
Given a measure-preserving action with a spectral gap of a connected simple Lie group $G$
on a (standard) probability space $(X,\mu)$ , for $ 1 < p < \infty$  there exist $c_p,\theta_p>0$ such that for the averaging operators  
$\pi_X(\beta_t)$ as in (\ref{eq:average0}), 
$$
\norm{\pi_X(\beta_t)f-\int_X fd\mu}_{p}\le C_p\,|B_t|^{-\theta_p}\norm{f}_p, \quad f\in L^p(X).
$$
Moreover, the parameters $c_p,\theta_p$ depend only on $p$ and the integrability exponent of the representation
$\pi_X|_{L^2_0(X)}$.
\end{thm}

It is a remarkable phenomenon that the rate of convergence  stated in Theorem \ref{spectral transfer}
depends only on the measure of the sets $B_t$.  Thus, $L^2$-convergence 
with a rate in a general ergodic action requires no geometric, regularity or stability assumptions of the family 
of sets in question, so that  the quantitative mean ergodic theorem is very robust. 

The pointwise ergodic theore, however,  is considerably more delicate, and for its validity some stability and
regularity assumptions on the family $B_t$ are indispensable. 
We turn now to discuss these issues in greater detail.

\begin{defi} \label{def:addmiss}
{\rm 
An increasing family of bounded measurable subset $B_t$, $t>0$, of $G$ is called 
{\it admissible} if 
there exists $c>0$ such that for all $t\ge t_0$  and $\vre\in (0,\vre_0)$,
\begin{align*}
{O}_\vre\cdot B_t\cdot {O}_\vre &\subset B_{t+c\vre},
\end{align*}
where $\mathcal{O}_\epsilon$ denotes the $\vre$-neighbourhood of identity with respect to a right
invariant metric on $G$, 
and
\begin{align*}
|B_{t+\vre}|&\le (1+c\,\vre)  |B_{t}|.
\end{align*}
}
\end{defi}

For instance, if $G$ is embedded as a closed subgroup of $\hbox{GL}_n(\mathbb{R})$ and 
$\|\cdot \|$ is norm on $\hbox{M}_n(\mathbb{R})$, 
one can show that a family 
$$
B_t=\{g\in G:\, \log\|g\|<t\}
$$
is admissible (see \cite[Appendix]{ems}, \cite[Ch.~7]{GN10}).

Under this mild regularity property, we also establish a quantitative pointwise ergodic theorem
\cite{MNS,GN10}:

\begin{thm}\label{th:pointwise}  
Consider a measure-preserving action with a spectral gap of a connected simple Lie group $G$
on a (standard) probability space, and let $\pi_X(\beta_t)$ be the averaging operators
with respect to an admissible family of sets.
Then  there exist $c_p,\theta_p>0$ such that for every $f\in L^p(X)$, $p>1$, and $\mu$-almost every $x\in X$,
$$
\left|\pi_X(\beta_t)f(x)-\int_X f\,d\mu\right|\le C_p(f,x)\,|B_t|^{-\theta'_p}, 
$$
where $\|C_p(f,\cdot)\|_p\le C_p'\, \|f\|_p$.
\end{thm}

One can also establish pointwise ergodic theorem for general actions
 without the spectral gap assumption, but in this case 
there is no estimate on the rate of convergence.
We refer to \cite{GN10} for a detailed account or ergodic theorem for semisimple groups.

\subsection{On spectral methods in ergodic theory} 
Before proceeding with the description of the applications of the mean and pointwise quantitative ergodic theorems stated above, let us make several very brief comments on some of the uses that spectral methods have found in ergodic theory. 
\subsubsection*{1.  Euclidean groups : singular averages} 
The methods of classical harmonic analysis on $\R^d$, including singular integrals, square function estimates, restriction theorems, fractional integrals and analytic interpolation have  been utilized  by several authors to prove ground-breaking pointwise ergodic theorems. Pioneering contribution here were Stein's ergodic theorem for power of a self adjoint operator \cite{St61}, and the celebrated spherical maximal theorem in $\R^d$, proved by Stein-Wainger for $d \ge 3$ \cite{SW78} and by Bourgain for $d=2$ \cite{Bo83}. The spherical maximal inequalities were used to prove a pointwise ergodic theorem for sphere averages, with the case $d \ge 3$ established by Jones \cite{Jo93}, and the case $d=2$ by Lacey \cite{La95}. Ergodic theorems for sphere averages on the Heisenberg group were developed in \cite{NT97}. 

For actions of the integers, the methods of discrete  harmonic analysis and especially estimates of exponential sums were developed to prove pointwise ergodic theorem for sparse sets of integers, such as square and primes,  in Bourgain's pioneering work on the subject \cite{Bo89}.  For actions of the integer lattices $\Z^d$, ergodic theorems were developed for discrete spheres and other averages by Magyar \cite{Mag02}. 
We remark that all the ergodic theorems just mentioned depend on the transference principle for the groups involved, namely they reduce the proof of the requisite maximal inequalities to the case of the action by convolution. 

\subsubsection*{2.  Free groups : radial averages}
The problem of establishing ergodic theorems for general non-amenable groups was raised already half a century ago by
Arnol'd and Krylov  \cite{AK63}.  An important case that has figured prominently in the theory ever since is that of the 
free (non-Abelian) group. Arnold and Krylov proved an equidistribution theorem for radial averages on dense free subgroups of isometries
of the unit sphere $\mathbb{S}^2$ via a spectral argument similar to Weyl's equidistribution theorem on
the circle.  
 Guivarc'h  has established a mean ergodic theorem for radial averages on the free group, using
 von-Neumann's original approach via the spectral theorem \cite{gu}. 
The pointwise ergodic theorem for general actions of the free group was proved in \cite{N1} for $L^2$-functions, 
and extended to function in $L^p$, $p > 1$, in \cite{NS1}. 
The distinctly non-amenable phenomenon of an ergodic theorem  with a quantitative estimate on the rate of convergence was realized by the  
celebrated Lubotzky-Phillips-Sarnak construction \cite{LPS1,LPS2} of a dense free group of isometries of $\mathbb{S}^2$ which has an optimal  spectral gap.

\subsubsection*{3.  Semisimple groups : ball and sphere averages}
An important and influential general spectral result on ergodic actions of semisimple Lie group has been
the Howe-Moore mixing theorem \cite{hm}, a consequence of the decay of matrix coefficients in unitary
representations of these groups.  
Tempelman (see \cite{Te2}) has used the Howe-Moore theorem to prove mean ergodic theorems for averages on
semisimple Lie group.

Some of the methods of classical harmonic analysis alluded to above were combined with unitary
representation theory of semisimple groups and used extensively in the development of pointwise ergodic
theorems for radial averages on semisimple Lie groups and some of their lattices. Pointwise ergodic
theorems for radial averages on the free groups were developed in \cite{N1,NS1}.  For  sphere and ball
averages on simple groups of real rank one pointwise ergodic theorems  were developed  in
\cite{N94b,N97,NS97}, and for complex groups in \cite{CN01}.   Pointwise ergodic theorems with a rate of
convergence to the ergodic mean were developed for radial averages on general semisimple Lie groups in
\cite{MNS} and for more general averages in \cite{N2}.

\section{Counting lattice points} \label{sec:counting}

Let $G$ be a locally compact second countable group,
and $\Gamma$ is a discrete subgroup with finite covolume.
We would be interested in determining the asymptotic behaviour of the number of elements
of $\Gamma$ contained in a family of increasing compact domains $B_t$ in $G$.
The problem of this type arises naturally in arithmetic geometry when ${\rm X}\subset \mathbb{C}^n$
is an affine algebraic variety, and one is interested in estimating the cardinality
of the set of integral points $x\in {\rm X}(\mathbb{Z})$ with bounded norm.
The techniques discussed in this section can be applied in particular  
when ${\rm G}\subset \hbox{GL}_n(\mathbb{C})$ is
an algebraic group defined over $\mathbb{Q}$, and $\Gamma={\rm G}(\mathbb{Z})$
is its arithmetic subgroup, and more generally to symmetric varieties of $\rm G$.

The simplest instance of the lattice counting problem is to estimate 
the number of integral vectors contained 
in a compact domain $B\subset \mathbb{R}^d$. In this case, it is easy to see that
$$
|\mathbb{Z}^d \cap B|=|B|+O(|\partial B|),
$$
and typically $\frac{|\partial B|}{|B|}\to 0$ as the volume $|B|\to \infty$, which leads
to an asymptotic formula. The analogous problem for more general groups presents
new significant challenges. In particular, for groups with exponential volume growth,
$|B|$ is comparable to $|\partial B|$, and the fundamental domain for $\Gamma$ in $G$
might be unbounded. 
Thus, in order to establish an asymptotic formula for $|\Gamma\cap B|$,
one is required to use more sophisticated analytic techniques.
A number of  methods have been developed to address this 
problem, which include, in particular:
\begin{itemize}
\item spectral expansion of the automorphic kernel
\cite{del, hu, se, pa, good, lp, mw92, drs, bmw, bgm},
\item decay of matrix coefficients \cite{mar, bar1, em, mau, bo},
\item symbolic coding and transfer operator techniques \cite{la, po},
\item the theory of unipotent flows \cite{ems}.
\end{itemize}
We refer to \cite{bab} for an extensive survey of the lattice counting problem.

Our approach to the lattice point counting problem in the domains $B_t$ is to study the averaging operators defined by $B_t$ in the action of $G$
on the homogeneous space $X=G/\Gamma$ equipped with the invariant probability measure $\mu$.
Given a family of domains $B_t$ in $G$, we consider
$$
\pi_X(\beta_t)f(x)=\frac{1}{|B_t|}\int_{B_t} f(g^{-1}x)\, dg,\quad f\in L^2(X).
$$
We now turn to  demonstrate that a (quantitative) mean ergodic theorem for these operators implies a 
(quantitative) solution of the lattice point counting problem for the famiy
of sets $B_t$, provided that it satisfies the regularity property of Definition \ref{def:addmiss}.

\begin{thm}\label{Thm:LP}
Let $G$ be a locally compact second countable group, $\Gamma$ a discrete lattice subgroup and assume that the family of sets $B_t$ is admissible. 
\begin{enumerate}
\item[(i)] Assume the averages $\pi_{G/\Gamma}(\beta_t)$ satisfy the mean ergodic theorem:
$$
\norm{\pi_X(\beta_t)f- \int_{G/\Gamma} f\, d\mu }_{2}\to 0,\quad f\in L^2(X).
$$ 
Then
$$
\abs{\Gamma\cap B_t}\sim 
\frac{|B_t|}{|G/\Gamma|}\quad\hbox{as}\quad t\to\infty.
$$
\item[(ii)] Assume the averages $\pi_{X}(\beta_t)$ satisfy the quantitative mean ergodic theorem:
\begin{equation}\label{eq:qqq}
\norm{\pi_{X}(\beta_t)f-\int_{X} f\, d\mu}_{2}\le C_2\, 
|B_t|^{-\theta }\norm{f}_{2},
\quad f\in L^2(X,\mu),
\end{equation}
with $C_2,\theta>0$.
Then for all sufficiently large $t$,\footnote{ By $\dim(G)$ here we mean the upper  local dimension, namely 
the infimum over $d>0$ such that $|{O}_\epsilon|\ge \epsilon^d$ for all sufficiently small
$\epsilon>0$,
where ${O}_\vare$ denotes the $\vare$-neighbourhood of identity.}
$$
|\Gamma\cap B_t|=\frac{|B_t|}{|G/\Gamma|}+O\left(|B_t|^{1-\frac{\theta}{\dim(G)+1}+\rho}\right),
$$ 
for every $\rho>0$,
where the implied constants depends only on the parameters in (\ref{eq:qqq}) and Definition
\ref{def:addmiss}, and the covolume of $\Gamma$.
\end{enumerate}
\end{thm}






%


%

Let us outline the proof of Theorem \ref{Thm:LP}.
We first observe that the quantity $|\Gamma\cap B_t|$
can be approximated by suitable averages on space $X$.
Recall that we denote by $O_\vare$ the $\vare$-neighbourhood of identity with respect 
to a right-invariant metric on $G$, which clearly satisfies $O_\vare^{-1}=O_\vare$.
We normalise the invariant measure on $G$ so that $\Gamma$ has covolume one.
Let $\chi_\vare$ be the characteristic function of 
this neighbourhood normalised to have integral one and
$f_\vare(g\Gamma)=\sum_{\gamma\in \Gamma} \chi_\vare(g\gamma)$.
Then $\int_X f\,d\mu=1$.
We claim that for $h\in {O}_\vare$,
\begin{equation}\label{eq:claim}
\int_{B_{t-c\vare}}f_\vare(g^{-1}h\Gamma)\,dg\le |\Gamma\cap B_t|\le 
 \int_{B_{t+c\vare}}f_\vare(g^{-1}h\Gamma)\,dg.
\end{equation}
Indeed, observe that 
$$
\int_{B_{t}}f_\vare(g^{-1}h\Gamma)\,dg= \sum_{\gamma\in\Gamma}
\int_{B_{t}}\chi_\vare(g^{-1}h\gamma)\,dg=\sum_{\gamma\in\Gamma}
\frac{|B_t\cap h\gamma {O}_\vare|}{|{O}_\vare|}.
$$
If $\gamma\in B_{t-c\epsilon}$, then by the admissibility property, 
$$
h\gamma {O}_\vare\subset O_\vare B_{t-c\vare}
O_\vare\subset B_t,
$$
so that
$$
|\Gamma\cap B_{t-c\vare}|
=\sum_{\gamma\in\Gamma\cap B_{t-c\vare}}
\frac{|B_t\cap h\gamma {O}_\vare|}{|{O}_\vare|}
\le 
\sum_{\gamma\in\Gamma}
\frac{|B_t\cap h\gamma {O}_\vare|}{|{O}_\vare|}.
$$
On the other hand, if $h\gamma {O}_\vare\cap B_t\ne \emptyset$, then
by the admissibility property,
$$
\gamma\in h^{-1}B_t{O}_\vare\subset O_\vare B_t O_\vare \subset B_{t+c\vare}.
$$
Hence, 
$$
\sum_{\gamma\in\Gamma}
\frac{|B_t\cap h\gamma{O}_\vare|}{|{O}_\vare|}
=
\sum_{\gamma\in\Gamma\cap B_{t+c\vare}}
\frac{|B_t\cap h\gamma {O}_\vare|}{|{O}_\vare|}
\le |\Gamma\cap B_{t+c\vare}|.
$$
This proves (\ref{eq:claim}).

It follows from the mean ergodic theorem that
for every $\delta>0$, 
$$
\mu(\{h\Gamma\in G/\Gamma:\, 
|\pi_{X}(\beta_t)f_\vare(h\Gamma) -1|>\delta\})
\rightarrow 0
$$
as $t\to\infty$. In particular, it follows for all sufficiently large $t$,
there exists $h_t\in {O}_\vare$ such that 
\begin{equation}\label{eq:aa}
|\pi_{X}(\beta_t)f_\vare(h_t\Gamma) -1|\le\delta.
\end{equation}
Combining this estimate with (\ref{eq:claim}), we deduce that
\begin{equation}\label{eq:bb}
|\Gamma\cap B_t|\le (1 + \delta)|B_{t+c\vare}|\le (1+\delta)(1+c^2\vare)|B_t|
\end{equation}
and similarly,
\begin{equation}\label{eq:cc}
|\Gamma\cap B_t|\ge  (1-\delta)(1-c^2\vare)|B_t|
\end{equation}
for sufficiently large $t$. Since these estimates hold for arbitrary $\delta,\vare>0$,
this proves the first part of Theorem \ref{Thm:LP}.

The proof of the quantitative estimate follows similar ideas, and is equally straightforward. 
We assume that $\vare$ is sufficiently small, and in particular the neighborhoods  $\mathcal{O}_\vare\gamma$,
$\gamma\in\Gamma$, are disjoint.
Then  since $|\mathcal{O}_\vare|\ge \vare^{n}$ with $n>\dim (G)$,
we deduce that  
$$
\norm{\chi_\vare}^2_2=|\mathcal{O}_\vare|^{-1}\le \vare^{-n}.
$$ 
The validity of \eqref{eq:qqq} now implies that 
$$
\mu(\{h\Gamma\in G/\Gamma:\, |\pi_{X}
(\beta_t)f_\vare(h\Gamma) -1|>\delta\})\le 
C_2^2\, \delta^{-2} \vare^{-n}  \abs{B_t}^{-2\theta }.
$$
Taking $\delta=2C_2\vare^{-n}|B_t|^{-\theta}$, we observe that the above measure
is strictly less than $\mu(\mathcal{O}_\vare\Gamma)=|\mathcal{O}_\vare|$,
 so that we can find $h_t\in \mathcal{O}_\vare$
satisfying \eqref{eq:aa} and deduce estimates \eqref{eq:bb}--\eqref{eq:cc} above.
Finally, to optimise the estimate, we take 
$\vare= |B_t|^{-\theta/(n+1)}$ with sufficiently large $t$. This leads to the second part of the theorem.

\vspace{0.2cm}

Let us note that Theorem \ref{Thm:LP} can be generalized significantly in two important respects. For full details we refer to \cite{gn_counting}.  
\begin{enumerate}
\item The quantitative lattice point count extends to  families of  sets $B_t$ much more general than admissible ones. 
It suffices that the families satisfy  the following weaker regularity property:
there exist $c,a>0$ such that for all $t \ge  t_0$ and  $\vare\in (0,\vare_0)$,
$$
\left|\bigcup_{u,v\in {O}_\vare}u B_tv\right|\le (1+c\,\vare^a)\left|\bigcap_{u,v\in {O}_\vare}
  uB_t v\right|,
$$
which we call \emph{H\"older well-roundedness}. 
This extends the solution for example to families of sector averages, families associated naturally 
with symmetric varieties, and many others. 
\item  The quantitative solution to the lattice point counting problem is uniform as the lattice varies in the space of lattice subgroups. All that is needed is that  the averaging operators satisfy a uniform norm decay estimate  in the space $L^2(G/\Gamma)$, and that $G/\Gamma$ contain an injective copy of a fixed neighborhood of the identity in $G$, as $\Gamma$ varies. In particular, uniformity holds over any set of finite index subgroups which satisfies property $(\tau)$.  Note that in the ergodic-theoretic approach we have taken, 
this feature of  uniformity is completely  obvious from the proof. 
\end{enumerate}

In particular,  combining the results of this section with Theorem \ref{spectral transfer} and  property $(\tau)$ for congruence subgroups, we deduce:

\begin{cor}\label{cor:c_count}
Let ${\rm G}\subset\hbox{\rm GL}_n(\CC)$ be a simply connected simple algebraic group defined over $\mathbb{Q}$
and $\Gamma(\ell)=\{\gamma\in {\rm G}(\mathbb{Z}):\, \gamma=\hbox{\rm id}\,\,\hbox{\rm mod}\,\ell\}$
be the family of congruence subgroups. Then for any H\"older well-rounded family of sets $B_t$ in 
${\rm G}(\mathbb{R})$, there exists $\delta>0$ such that for sufficiently large $t$,
$$
|\Gamma(\ell)\cap B_t|=\frac{|B_t|}{|{\rm G}(\mathbb{R})/\Gamma(\ell)|}+O(|B_t|^{1-\delta}),
$$
where the implied constant is independent of $\ell$. 
\end{cor}

This estimate will play an indispensable role in the next section.

\section{Almost prime points on homogeneous algebraic varieties}\label{sec:prime}
In the present section we will describe a solution to the problem of sifting the integral points on certain homogeneous algebraic varieties to obtain almost prime points,  following \cite{ns, gn_lifting}.  As will be explained in due course, the uniformity in the lattice point counting problem established  in the previous section will play a crucial role.  But let's begin with some remarks on the origins of this problem. 

\subsection{Sifting along finite-index subgroups}
The prime number theorem  describes the asymptotic distribution of the set $\mathbb{P}=\{2,3,5,7,11,\ldots\}$
of prime numbers:
\begin{equation}\label{eq:prime}
|\{p\in\mathbb{P}:\, p\le t\}|\sim \frac{t}{\log t}\quad\hbox{as $t\to \infty$.}
\end{equation}
More generally, let us consider an integral irreducible polynomial $f$.
We assume that the leading coefficient of $f$ is positive and $\gcd(f(\mathbb{Z}))=1$.
Is it true that $f(x)$ is prime for infinitely many integers $x$?
The question is answered for linear polynomials by the Dirichlet theorem about primes
in arithmetic progressions, but it turns out to be extremely difficult
for polynomials of higher degrees. 
A natural approximation to this problem is the question whether 
for a fixed $r\ge 1$, the values $f(x)$
belongs to the set  $\mathbb{P}_r$ of almost primes infinitely often, where 
$\mathbb{P}_r$ denotes the set of integers having at most $r$ prime factors.
This problem has been extensively studied using sieve methods.
For example, it has been shown (see \cite[Sec.~9.5]{hr})
that when $r= \deg(f)+1$, there is a lower estimate of correct order:
\begin{equation}\label{eq:almost prime}
|\{x\in \mathbb{Z}:\,\, f(x)\in \mathbb{P}_r, |x|\le t\}|\ge \hbox{const}\cdot\frac{t}{\log t}.
\end{equation}
Although it is widely believed that under the above necessary conditions on $f$,
the same estimate should holds for the set of primes as well, this is
currently out of reach. In particular, it is not known weither there are
infinitely many primes of the form $n^2+1$ (see \cite{i} for the best result in this direction).

A far-reaching conjectural generalisation of this problem was proposed by Bourgain, Gamburd and Sarnak (see  \cite{s_primes} and \cite{bgs2}), and termed the saturation problem, as follows. Given  
\begin{itemize}
\item a Zariski dense subgroup $\Gamma$ of $\hbox{GL}_d(\mathbb{\Z})$,
\item a linear representation $\rho:\hbox{GL}_d(\CC)\to \hbox{GL}_n(\CC)$ defined over $\mathbb{Q}$
such that $\rho(\Gamma)\subset \hbox{GL}_n(\mathbb{Z})$,
\item a vector $v\in \mathbb{Z}^n$, 
\item a polynomial $f\in\mathbb{Z}[x_1,\ldots,x_n]$,
\end{itemize}
is it true that $f(x)$, with $x\in\rho(\Gamma) v$, is $r$-prime infinitely often, 
and more generally, is the set of such points $x$ Zariski dense in the Zariski closure of $\rho(\Gamma) v$?

One can also ask for a quantitative version of this question. Setting
$$
\mathcal{O}(t)=\{x\in \Gamma v:\, \|x\|\le t\},
$$
is it true that 
\begin{equation}\label{eq:almost_prime}
|\{x\in \mathcal{O}(t):\, f(x)\in\mathbb{P}_r \}|\ge \hbox{const} \cdot \frac{|\mathcal{O}(t)|}{\log t}\quad \hbox{as $t\to\infty$,}
\end{equation}
in analogy with  \eqref{eq:almost prime}?

It was a profound idea of P.~Sarnak that classical sieve methods can be applied in this
noncommutative setting provided that one is able to establish an asymptotic estimate 
for the cardinalities of the sets
\begin{equation}\label{eq:cong0}
\mathcal{O}(t,\ell)=\{x\in \Gamma v:\,\, \|x\|\le t,\, f(x)=0\,\hbox{mod}\, \ell\}
\end{equation}
as $t\to\infty$.
 It is crucial for the sieve method to succeed that this estimate is sufficiently uniform over $\ell$. Letting 
\begin{equation}\label{eq:congruence0}
\Gamma(\ell)=\{\gamma\in \Gamma:\, \gamma=\hbox{id}\,\hbox{mod}\, \ell\}
\end{equation}
denote the congruence subgroup of level $\ell$, it is clear that 
 $f(\gamma\Gamma(\ell)v)=f(\gamma v)\,\hbox{mod}\, \ell$ for  $\gamma\in \Gamma/\Gamma(\ell)$.
 Thus in order to estimate $|\mathcal{O}(t,\ell)|$, it is sufficient to estimate
cardinalities of the sets
\begin{equation}\label{eq:cong}
\{x\in \gamma\Gamma(\ell) v:\, \|x\|\le t\}\quad\text{ for } \gamma\in \Gamma/\Gamma(\ell)\,.
\end{equation}

Let us now consider the important special case of lattice actions on principal homogeneous spaces, where $\mathcal{O}=\Gamma v\subset Gv$, $\Gamma$ is a lattice in $G$, and the stability group of $v$ in $G$ is finite. Then estimating (\ref{eq:cong}) amounts to nothing other than the lattice point counting problem, specifically in the intersection of norm balls  with  translates of the orbits of congruence subgroups of $\Gamma$.  It is crucial for the sieve method to succeed that  the solution to the lattice point counting problem in such orbits be uniform over the congruence subgroups. The uniform counting property holds in these sets  when  $\Gamma$ is an arithmetic
subgroup of a semisimple algebraic group $\rm G$ (for instance, when $\Gamma=\hbox{SL}_d(\mathbb{Z})$) \cite{ns}. 
Thus  the required estimates for the sieve method is ultimately derived from the uniform spectral gap
property, namely property $(\tau)$ discussed in Section \ref{sec:ug}.
Indeed, property $(\tau)$ gives a uniform spectral estimate for the ball averages acting  on the spaces $L^2({\rm G}(\mathbb{R})/\Gamma(\ell))$, 
which ultimately implies a uniform estimate in the lattice point counting problem.

Establishing a Zariski dense subset of almost prime points on homogeneous varieties  requires the more general uniform lattice point counting results of Corollary \ref{cor:c_count}. Before stating such a result, 
 let us illustrate how the estimates on $|\mathcal{O}(t,\ell)|$ comes into play
on a simple sieve based on the Mobius function.
Recall that the Mobius function $\mu(\ell)$ is equal to $(-1)^r$, if $\ell$ is a product of $r$ distinct prime
factors, and is $0$, otherwise, and the Mobius inversion formula
$$
\sum_{\ell:\, \ell|n} \mu(\ell)=\left\{
\begin{tabular}{l}
$1$ if $n=1$,\\
$0$ if $n>1$.
\end{tabular} \right.
$$
Let $P_z$ be the product of all prime numbers up to $z$ and $S(t,z)$ denotes the cardinality of
the set of $x\in \mathcal{O}(t)$ such that $f(x)$ is coprime to $P_z$. Then
\begin{align}\label{eq:sieving}
S(t,z)&=\sum_{k: \gcd(k,P_z)=1} |\mathcal{O}(t)\cap \{f=k\}|\\ \nonumber
&=\sum_k \left(\sum_{\ell:\, \ell|\gcd(k,P_z)} \mu(\ell)\right)|\mathcal{O}(t)\cap \{f=k\}|\\ \nonumber
&=\sum_{\ell:\, \ell| P_z} \mu(\ell) \sum_{k:\, \ell|k} |\mathcal{O}(t)\cap \{f=k\}|\\ \nonumber
&=\sum_{\ell:\, \ell| P_z} \mu(\ell) |\mathcal{O}(t,\ell)|.
\end{align}
Now suppose that for some $\kappa,\delta>0$, we have an estimate
\begin{equation}\label{eq:main_est}
|\mathcal{O}(t,\ell)|=\rho(\ell)|\mathcal{O}(t)|+
O\left(\ell^\kappa |\mathcal{O}(t)|^{1-\delta}\right)
\end{equation}
with some uniform $\kappa,\delta>0$, where
$$
\rho(\ell)=\frac{|\{\Gamma v\,\hbox{mod}\, \ell\}\cap \{f=0\,\hbox{mod}\, \ell\}|}{|\{\Gamma v\,\hbox{mod}\, \ell\}|}.
$$
When $z$ is sufficiently small compared with $t$, the above computation leads to
favourable estimate for $S(t,z)$. However, to derive (\ref{eq:almost_prime}),
one needs to estimate $S(t,z)$ when $z=t^\alpha$ with some $\alpha>0$, and the above sieve
with the Mobius weights is not sufficiently efficient for our purpose.
Fortunately, starting with the works of V.~Brunn in 1920's more efficient sieving weights
have been developed. We refer to recent surveys \cite{f,mo,gr,mol} and monographs
\cite{fi,g,h,hr,ik} for a comprehensive discussion of sieve methods.
These techniques allow to deduce a favourable estimate on $S(t,z)$ with $z=t^\alpha$ with $\alpha>0$
provided that one knows (\ref{eq:main_est}) and asymptotic estimates of the
coefficients $\rho(\ell)$. The latter estimates can be deduced from the strong
approximation property and results on the number of points on varieties over
finite fields. Ultimately the main difficulty lies in establishing
(\ref{eq:main_est}), and the crucial initial input is the uniformity in the lattice point counting which follows from property $(\tau)$.

More generally, it is possible to apply this method to the case of orbits on symmetric varieties of $G$. In particular, let us mention the following result \cite[Example~1.11]{gn_lifting}, 
which shows that quadratic surfaces contain many integral points with almost prime coordinates
and implies that the set of almost prime points is not contained in any proper algebraic subset.
This example is a partial case of Theorem \ref{th:primes} below.
Let 
$Q(x)=\sum_{i,j=1}^d a_{ij}x_ix_j$
be a nondegenerate indefinite integral quadratic form in $d$ variables.
We assume that $d\ge 4$ (see \cite{ls} for $d=3$ case)
and denote by $\hbox{Spin}_Q$ the spinor group of $Q$.
We fix $v\in \mathbb{Z}^d$ such that $n=Q(v)\ne 0$ and consider the orbit $\mathcal{O}=\hbox{Spin}_{Q}(\mathbb{Z})v$
which lies in the quadratic surface ${\rm X}=\{Q(x)=n\}$.
Assuming that the polynomial function $f$ is absolutely irreducible on ${\rm X}$
and $\gcd(f(x): x\in \mathcal{O})=1$,
there exists explicit $r=r(d,\deg(f))\ge 1$ such that (\ref{eq:almost_prime}) holds.

Turning to sifting on general symmetric varieties, recall that Corollary \ref{cor:c_count} 
establishes uniform asymptotics for 
the number of elements of $\Gamma(\ell)$ contained 
in a family of H\"older well-rounded subsets $B_t$ of ${\rm G}(\mathbb{R})$.
More generally, as noted already in the discussion of almost prime points above,  one can also deal of cosets of $\Gamma(\ell)$ and establish
the following estimate. For every $\gamma\in\Gamma/\Gamma(\ell)$,
\begin{equation}\label{eq:count_cong}
|\gamma\Gamma(\ell)\cap B_t|=\frac{|B_t|}{|{\rm G}(\mathbb{R})/\Gamma(\ell)|}
+O\left(|B_t|^{1-\delta}\right),
\end{equation}
where $\delta>0$ depends only 
on the supremum of the integrability exponents in the representations $L^2_0(G/\Gamma(\ell))$, the dimension of $\rm G$,
and the H\"older exponent of $B_t$'s, and the implicit constant is independent of $\ell$ and $\gamma\in\Gamma/\Gamma(\ell)$.
Estimate (\ref{eq:count_cong}) can be used directly for sieving for almost primes
on the arithmetic groups $\Gamma$ itself, or more generally in its orbits in principal homogeneous spaces, and this  was carried out in \cite{ns}.
To sieve for almost primes in orbits $\mathcal{O}=\Gamma v$, one may attempt 
to construct a H\"older well-rounded family $\{B_t\}$ of sets in ${\rm G}(\mathbb{R})$
such that the map
$$
\Gamma\cap B_t \to \mathcal{O}(t):\gamma\mapsto \gamma v
$$
is onto and has fibers with uniformly bounded cardinalities. Then using (\ref{eq:count_cong}),
one would show that there exists $r\ge 1$ such that
$$
|\{\gamma\in \Gamma\cap B_t:\, f(\gamma v)\in\mathbb{P}_r \}|
\ge \hbox{const} \cdot \frac{|\Gamma\cap B_t|}{\log t}\quad \hbox{as $t\to\infty$,}
$$
which implies (\ref{eq:almost_prime}).
This strategy has been realised in \cite{gn_lifting} when the stabiliser of $v$ in $\rm G$ is 
a symmetric subgroup, namely the set of fixed
points of an involution. It leads to the following  result (see \cite[Theorem~1.10]{gn_lifting}) on the saturation problem for symmetric varieties.

\begin{thm}\label{th:primes}
Let ${\rm G}\subset\hbox{\rm GL}_n(\CC)$ be a connected $\mathbb{Q}$-simple simply connected
algebraic group, $\Gamma={\rm G}(\mathbb{Z})$, 
and $v\in\mathbb{Z}^n$ be such that $\hbox{\rm Stab}_{\rm G}(v)$ is connected and symmetric.
Then given any polynomial $f:\mathbb{Z}^n\to \mathbb{Z}$ such
the function $g\mapsto f(gv)$ is absolutely irreducible on $\rm G$ and $\gcd(f(x):x\in\mathcal{O})=1$,
there exists explicit $r=r({\rm G},v,\deg(f))\ge 1$ for which
the estimate (\ref{eq:almost_prime}) holds.
\end{thm}

\subsection{Bounds  on integral point satisfying congruence conditions }

Another fundamental question regarding the distribution of prime numbers is the Linnik problem.
By Dirichlet's theorem, for every coprime $b,\ell\in \mathbb{N}$ one can find a prime number $p$
satisfying 
$$
p=b\,\,\hbox{mod}\, \ell.
$$
 Y.~Linnik \cite{l1,l2} raised the question whether this congruence can be solved
effectively, and showed that there exist $\ell_0,\sigma>0$ such that one can find
a prime $p$ satisfying $p\le \ell^\sigma$ for all $\ell\ge \ell_0$.
This problem admits a sweeping generalization to the  noncommutative set-up,  namely estimating the norm
of  a minimal prime solution to a system of polynomial equations  satisfying a congruence condition. The latter problem is out of reach at this time,  but it is possible to give a general solution to the problem of bounding  the norm of an almost prime solution, as follows. 

Let $\mathcal{O}\subset\mathbb{Z}^d$ be an orbit of an arithmetic group 
and $f:\mathcal{O}\to \mathbb{Z}$ a polynomial function.
We would like to find a small solution of the congruence
\begin{equation}\label{eq:congruence}
f(x)=b\,\,\hbox{mod}\,\ell,\quad x\in\mathcal{O},
\end{equation}
with $f(x)$ being almost prime.
Using the standard sieving techniques combined with the estimate (\ref{eq:count_cong}),
we prove the following 
(see \cite[Theorem~1.13]{gn_lifting}):

\begin{thm}\label{th:linnik}
Let ${\rm G}\subset\hbox{\rm GL}_n(\CC)$ be a connected $\mathbb{Q}$-simple simply connected
algebraic group, $\Gamma={\rm G}(\mathbb{Z})$,  and $v\in\mathbb{Z}^n$.
Then given a polynomial $f:\mathbb{Z}^n\to \mathbb{Z}$ such
the function $g\mapsto f(gv)$ is absolutely irreducible on $\rm G$ and $\gcd(f(x):x\in\mathcal{O})=1$,
there exist explicit $r,\ell_0,\sigma>0$ such that for every coprime $b,\ell\in\mathbb{N}$
satisfying $\ell\ge \ell_0$ and $b\in f(\mathcal{O})\,\hbox{mod}\, \ell$, one can find 
$x$ with $\|x\|\le \ell^\sigma$ satisfying (\ref{eq:congruence}) such that $f(x)$ is $r$-prime.
\end{thm}

The subject of sieving along orbits of noncomutative groups
is a rapidly developing area now. In particular, we mentionthe extensive  developments regarding Apollonian packings
\cite{s_ap0,s_ap,k_o2,bf}
and orbits of ``thin'' subgroups more generally \cite{bgs1,ko,bgs2,bk,k_o1}.
We refer to \cite{kont,oh} for recent comprehensive surveys.

\section{Intrinsic Diophantine approximation on homogeneous varieties}\label{sec:dioph}
\subsection{Diophantine approximation and ergodic theory} 
In the theory of Diophantine approximation one is interested in finding
an efficient rational approximation for vectors $x\in \mathbb{R}^d$.
More explicitly, one would like to know for what range of parameters $\epsilon$
and $R$, the system of inequalities
$$
\|x-r\|\le \epsilon\quad\hbox{and}\quad \hbox{D}(r)\le R
$$
has a solution $r\in\mathbb{Q}^d$, where $\|\cdot\|$ denotes the maximum norm, and $\hbox{D}(r)$
denotes the denominator of a rational number $r$ written in reduced form.
In this section we consider a more general question regarding intrinsic Diophantine approximation on
algebraic sets, namely approximation by rational points on the variety itself. Let 
$$
{\rm X}=\{f_1(x)=\cdots =f_k (x)=0\}\subset \mathbb{C}^n
$$
be an algebraic set  defined by a family of polynomial with rational coefficients.
For a real points $x$ in ${\rm X}(\mathbb{R})$ which is contained in the closure of 
the set ${\rm X}(\mathbb{Q})$ of rational points, we would like to produce
an estimate that quantifies density. This amounts to solving the system of inequalities 
\begin{equation}\label{eq:dioph}
\|x-r\|\le \epsilon\quad\hbox{and}\quad \hbox{D}(r)\le \epsilon^{-\kappa}
\end{equation}
with $r\in {\rm X}(\mathbb{Q})$.
The problem of Diophantine approximation on general algebraic varieties by rational points on it was raised in the 1965 survey by
Lang \cite{L65}.
This question cannot be addressed in full generality
because it is not even known how to determine
whether ${\rm X}(\mathbb{Q})$ is finite or infinite.
Here we propose an approach to it when ${\rm X}(\mathbb{Q})$ has an additional structure, namely when ${\rm X}(\mathbb{Q})$ is equipped with a group action.
We note that when ${\rm X}$ is an elliptic curve 
(or, more generally, an abelian variety) the problem of Diophantine approximation on $\rm X$
has been studied by M. Waldschmidt in \cite{w}. We 
have considered the problem of Diophantine approximation on an algebraic set $\rm X$
equipped with a transitive action of a semisimple algebraic group $\rm G$ defined over $\mathbb{Q}$.
A key new feature of our approach is that it allows approximation by rational points satisfying an arbitrary set of integrality constraints, and not only by the set of all rational points. 
To simplify our notation for this exposition, we consider
the problem of Diophantine approximation by the rational points in ${\rm X}(\mathbb{Z}[1/p])$
where $p$ is prime. The problem of Diophantine approximation by points in ${\rm X}(\mathbb{Z}[1/p_1,\dots, 1/p_s])$  (say), or by ${\rm X}(\mathbb{Q})$
can be handled similarly, but uses analysis on adelic spaces 
which would require much more elaborate notation.

Let us illustrate our general result stated below by the following example, mentioned in the introduction.  
Let ${\rm X}$ be a rational two-dimensional ellipsoid such that ${\rm X}(\mathbb{Z}[1/p])$ is not discrete
in ${\rm X}(\mathbb{R})$.
Then we prove that 
\begin{itemize}
\item for almost every $x\in {\rm X}(\mathbb{R})$, every $\kappa>2$, and $\epsilon\in (0,\epsilon_0(x,\kappa))$,
the system of inequalities
\begin{equation}\label{eq:sphere1}
\|x-r\|\le \epsilon\quad\hbox{and}\quad \hbox{D}(r)\le \epsilon^{-\kappa}
\end{equation}
has a solution $r\in {\rm X}(\mathbb{Z}[1/p])$.
\item for every $x\in {\rm X}(\mathbb{R})$, every $\kappa>4$, and $\epsilon\in (0,\epsilon_0(\kappa))$,
the system of inequalities
\begin{equation}\label{eq:sphere2}
\|x-r\|\le \epsilon\quad\hbox{and}\quad \hbox{D}(r)\le \epsilon^{-\kappa}
\end{equation}
has a solution $r\in {\rm X}(\mathbb{Z}[1/p])$.
\end{itemize}
We note that the estimate (\ref{eq:sphere1})
is essentially the best possible.
Indeed, the number of $r\in {\rm X}(\mathbb{Z}[1/p])$ with $\hbox{D}(r)\le R$
grows as $\hbox{const}\cdot R$ as $R\to\infty$ and $\dim ({\rm X})=2$,
so that one can deduce from the pigeon-hole principle that 
$\kappa\ge 2$ in (\ref{eq:sphere1}). 

For a general algebraic set ${\rm X}$, we set
$$
a_p({\rm X})=\sup_{\hbox{\tiny compact $\Omega\subset {\rm X}(\mathbb{R})$}} \limsup_{R\to\infty}
\frac{\log |\{r\in \Omega\cap {\rm X}(\mathbb{Z}[1/p]):\, \hbox{D}(r)\le R \}|}{\log(R)}.
$$
Then again it follows from the pigeon-hole that the exponent $\kappa=\frac{\dim({\rm X})}{a_p({\rm X})}$
in (\ref{eq:dioph}) is the best possible.

We note that rational ellipsoids constitute the only example for which we are aware of previous results  in the literature, as follows. 
Diophantine approximation by points in ${\rm X}(\Z[\frac1p])$ for ${\rm X}=S^2$ and ${\rm X}=S^3$ can be deduced from the celebrated construction by 
Lubotzky-Phillips-Sarnak of a dense subgroup of $\hbox{SO}(3,\R)$ with entries in $\Z[\frac{1}{p}]$ possessing the optimal spectral gap \cite{LPS1,LPS2}. 
A rate of uniform approximation by all rational points in the spheres $S^d$, $d\ge 2$ was established by Schmutz \cite{Sch08}, using elementary methods based on rational parametrizations of the spheres, the rate being given by $\kappa=2\log_2(d+1)$. Duke \cite{Du03} has established equidistribution of the set of all rational points on $S^2$, from which a rate of uniform approximation can be derived, and the same method applies to spheres of any dimension. Finally, we mention the remarkable recent results of Kleinbock and Merrill 
\cite{KlMe13}, which give the best possible results for approximation by all rational points on spheres of any dimension, establishing complete analogs of Dirichlet's and Khinchin's theorems in this case. 

Turning now to a general variety, in order to study the density of the set ${\rm X}(\mathbb{Z}[1/p])$ in ${\rm X}(\mathbb{R})$, we construct
a dynamical system which encodes it in some sense. This turns out possible
because the set ${\rm X}(\mathbb{Z}[1/p])$ is parametrised by the orbits of the group ${\rm G}(\mathbb{Z}[1/p])$.
It is crucial for the classical theory of Diophantine approximation
that $\mathbb{Z}^d$ forms a lattice $\mathbb{R}^d$, and in order to have a similar framework
in our setting we consider the topological  group ${\rm G}(\mathbb{R})\times {\rm G}(\mathbb{Q}_p)$ 
(where $\mathbb{Q}_p$ denotes the field of $p$-adic number) equipped with an invariant measure. 
Then the group ${\rm G}(\mathbb{Z}[1/p])$ embeds diagonally in
${\rm G}(\mathbb{R})\times {\rm G}(\mathbb{Q}_p)$ as a discrete subgroup with finite covolume.
We consider the space 
$$
Y=({\rm G}(\mathbb{R})\times {\rm G}(\mathbb{Q}_p))/{\rm G}(\mathbb{Z}[1/p]).
$$
The group ${\rm G}(\mathbb{Q}_p)$ acts naturally on this space preserving the finite measure.
This is the dynamical system that plays a crucial role in our considerations.
It turns out that from the dynamical systems point of view, the original problem of Diophantine approximation
corresponds to the {\it shrinking target property}.
Suppose that we have an increasing sequence of compact subsets $B_n$ of ${\rm G}(\mathbb{Q}_p)$ and
a shrinking sequence of neighbourhoods $O_{\vre}$ in $Y$. The shrinking target
problem asks for what range of parameters $n$ and $\epsilon$ the orbits $B_n y$ with $y\in Y$
reach the neighbourhoods ${O}_\epsilon$.
While in the classical theory of dynamical systems, this question is asked for a 
one-parameter flow, here we are required to answer a similar question for orbits of a large groups
${\rm G}(\mathbb{Q}_p)$. 
The connection between 
Diophantine approximation in $\mathbb{R}^d$ and the shrinking target property of one-parameter flows
have been previously exploited by S. G. Dani, D.~Kleinbock, and G.~Margulis \cite{d1,d,d2,kl,kl1,km1,kl2,km2}.
These developments are surveyed in \cite{kl3,mar,kl4}.
Also in the context of negatively curved manifolds, this connection  
was explored by S. Hersonsky and F. Paulin \cite{hp1,hp2,hp3}.

\begin{figure}[h]
\includegraphics[height=2in]{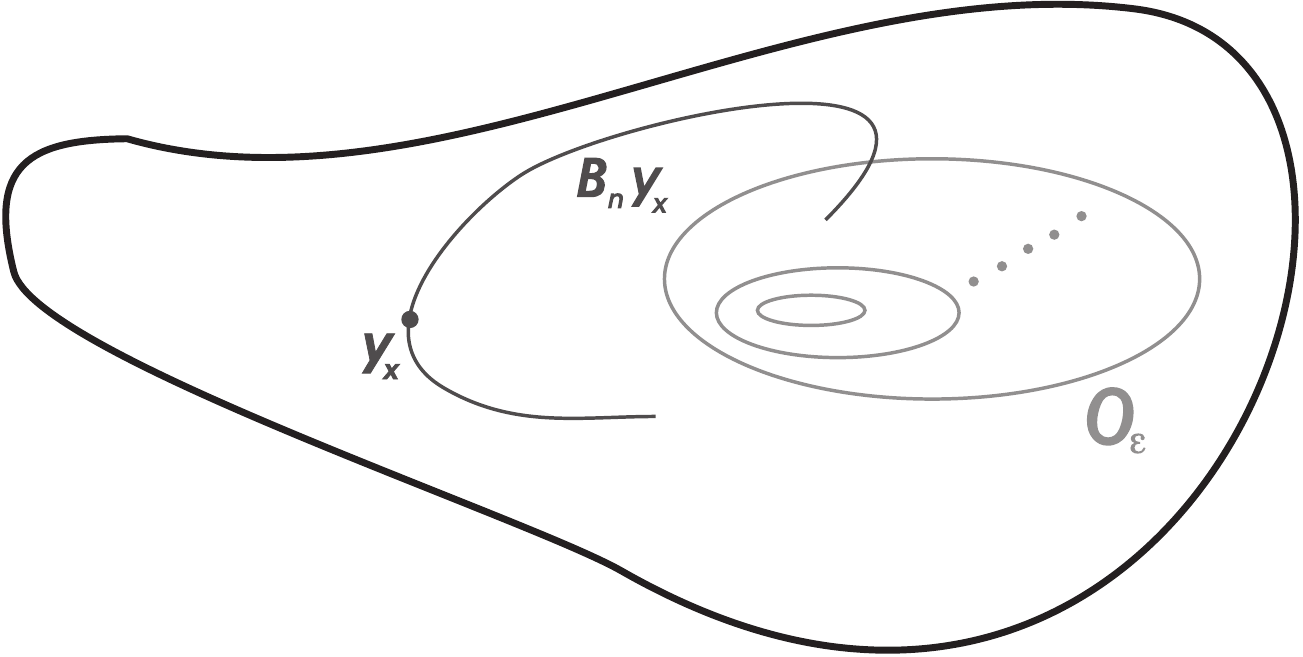}
\caption{Shrinking target property}
\end{figure}

\subsection{Diophantine approximation and the duality principle}
To address the problem of Diophantine approximation in our setting
we take 
\begin{equation}\label{eq:b_n}
B_n=\{g\in G(\mathbb{Q}_p):\, \|g\|_p\le p^n\}
\end{equation}
and construct a shrinking family of neighbourhoods $O_\vre$ in $Y$.
An important property of  $O_\vre$ is that it models the family of neighbourhoods $\|x-r\|\le \epsilon$ in ${\rm X}(\mathbb{R})$
and in particular $|O_\epsilon|$ decays as $\hbox{const}\cdot\epsilon^{\dim(X)}$
when $\vre\to 0^+$.
We show that every point $x\in X(\mathbb{R})$ is naturally associated  a point $y_x\in Y$, and 
the system of inequalities
$$
\|x-r\|\le \epsilon,\quad \hbox{D}(r)\le p^n
$$
has a solution $r\in {\rm X}(\mathbb{Z}[1/p])$ provided that
\begin{equation}\label{eq:shrink}
B_n^{-1}\, y_x\cap O_\epsilon\ne \emptyset.
\end{equation}
This leads to a useful connection 
(summarised in Table 1) between the 
problem of Diophantine approximation on $\rm X$
and the distribution of orbits of $Y$.
In the classical setting an analogous correspondence was observed by S. G. Dani in \cite{d1}
and called the Dani correspondence in \cite{km1}.
Since our construction is somewhat similar, we call it the {\it generalised Dani correspondence}, see Table 1 below.

\begin{table}[h]
\begin{tabular}{|c|c|}
\hline
&\\
Diophantine approximation & Dynamics \\
&\\
\hline\hline
&\\
${\rm X}(\mathbb{Z}[1/p])\subset {\rm X}(\mathbb{R})$ &
${\rm G}(\mathbb{Q}_p) \curvearrowright Y $ \\
&\\
\hline
& \\
$x\in {\rm X}(\mathbb{R})$ & $y_x\in Y$ \\
&\\
\hline
& \\
$\|x-r\|\le \epsilon$ & $O_\epsilon\subset Y$\\
&\\
\hline
&\\
$\left\{
\begin{tabular}{l}
$\|x-r\|\le \epsilon$,\\
$\hbox{den}(r)\le p^n$
\end{tabular}\right.$ has a $\mathbb{Z}[1/p]$-solution & $B_n^{-1}\, y_x \cap O_\epsilon\ne \emptyset$\\
&\\
\hline
\end{tabular}
\vspace{0.5cm}
\caption{Generalised Dani correspondence}
\end{table}

\subsection{Spectral estimates}

We prove the shrinking target property by analyzing the behavior of suitable averaging operators
on the space $Y$.
Let $\beta_n$ denote the uniform probability measure supported on the set $B_n\subset {\rm
  G}(\mathbb{Q}_p)$. Then we have the sequence of averaging operators 
$$
\pi_Y(\beta_n): L^2(Y)\to L^2(Y)
$$ 
defined by
$$
\pi_Y(\beta_n)\phi(y)=\int_{{\rm G}(\mathbb{Q}_p)} \phi(g^{-1}y)\,d\beta_n(g),\quad \phi\in L^2(Y).
$$
Ideally, one hopes that as in the classical von Neumann ergodic theorem, the averages
$\pi_Y(\beta_n)\phi$ converge to the space average $\int_Y \phi$ as $n\to\infty$, but this is not always true
due to presence of nontrivial unitary characters of ${\rm G}(\mathbb{Q}_p)$ in the space $L^2(Y)$.
We denote by $L^2_{00}(Y)$ the subspace of $L^2(Y)$ orthogonal to all characters of 
${\rm G}(\mathbb{R})\times{\rm G}(\mathbb{Q}_p)$. 
Then one can show that $\pi_Y(\beta_n)\phi\to 0$ in $L^2$-norm for every $\phi\in
L^2_{00}(Y)$. The rate of convergence here will play a crucial role.
This rate is deduce from a variant the spectral gap property, which was discussed in Section  4 and Section \ref{sec:prime},
and is determined by a suitable integrability exponent.
We fix a good maximal compact subgroup $U_p$ of ${\rm G}(\mathbb{Q}_p)$.
For functions $\phi,\psi\in L^2(Y)$, we denote
by $c_{\phi,\psi}(g)=\int_Y \phi(g y)\overline{\psi(y)}\, dy$ the corresponding matrix
coefficient.
Similarly to (\ref{eq:integ}) we define the {\it spherical integrability exponent} 
$$
q_p({\rm G})=\inf\left\{ q>0: 
\, \hbox{$c_{\phi,\psi}\in L^q({\rm G}(\mathbb{Q}_p))$
for all $U_p$-inv. $\phi,\psi\in L^2_{00}(Y)$}
\right\}.
$$
It follows from the result of L.~Clozel \cite{clo} that the integrability exponent
$q_p({\rm G})$ is finite, and we establish the following quantitative mean ergodic theorem:
for every sequence $B_n$ of $U_p$-bi-invariant compact subsets of ${\rm G}(\mathbb{Q}_p)$, 
every $\phi\in L^2_{00}(Y)$ and $\delta>0$,
\begin{equation}\label{eq:mean}
\|\pi_Y(\beta_n)(\phi)\|_2\le \hbox{const}(\delta)\cdot |B_n|^{-\frac{1}{q_p({\rm G})}+\delta}\|\phi\|_2.
\end{equation}
This mean ergodic theorem is vital for establishing the required shrinking target property.
We note that in order to derive optimal (or close to optimal) exponents for Diophantine
approximation in (\ref{eq:dioph}) it is crucial that the error term is controlled by
the $L^2$-norm rather than a Sobolev norm, which compromises the rate.
The actual value of the exponent of $|B_n|$ in (\ref{eq:mean}) is the crucial parameter which ultimately controls the quality  of  Diophantine approximation, 
as elaborated further below. 
Thus this approach provides a connection between the generalised Ramanujan conjectures
and the problem of Diophantine approximation on homogeneous varieties.

Now coming back to the original Diophantine approximation problem,
we apply estimate (\ref{eq:mean}) to the averages supported on the sets $B_n$ as in (\ref{eq:b_n})
(or, more precisely, on the $U_p$-bi-invariant sets $U_pB_nU_p$). 
We also estimate the behaviour of these averages on the orthogonal complement
of $L^2_{00}(Y)$ in $L^2(Y)$. Since the orthogonal complement 
is finite-dimensional, 
this does not present a particular challenge.
Assuming that
\begin{equation*}
|O_\epsilon|\ge\hbox{const}\cdot |B_n|^{-\frac{2}{q_p(G)}+\delta}
\end{equation*}
with fixed $\delta>0$, we use (\ref{eq:mean}) to produce an upper estimate
on the measure of the set $y\in Y$ for which $B_n^{-1} y\cap O_\vre=\emptyset$.
Then using a Borel--Cantelli-type argument, we deduce a shrinking target property
that holds for almost every $y\in Y$. This implies a Diophantine approximation result
via the generalised Dani correspondence (Table 1). If we impose a stronger assumption
that 
\begin{equation*}
|O_\epsilon|\ge\hbox{const}\cdot |B_n|^{-\frac{1}{q_p(G)}+\delta}
\end{equation*}
with some fixed $\delta>0$, then the above argument can be modified to
show that $B_n^{-1} y\cap O_\vre=\emptyset$ for all $y\in Y$ provided that
$\vre$ is sufficiently small. This produces a Diophantine approximation result
which is valid for all points in $\overline{{\rm X}(\mathbb{Z}[1/p])}$.

This completes the outline of the following result, which is part of a more general theorem stated in \cite[Theorem~1.6]{ggn}: 

\begin{thm}\label{th:dioph}
Let $\rm X$ be an algebraic set defined over $\mathbb{Q}$ equipped with a transitive action
of connected almost simple algebraic group $\rm G$ defined over $\mathbb{Q}$. We assume that
${\rm X}(\mathbb{Z}[1/p])$ is not discrete in ${\rm X}(\mathbb{R})$\footnote{One can show that under these assumtions
$\overline{{\rm X}(\mathbb{Z}[1/p])}$ is open in ${\rm X}(\mathbb{R})$.}.
Then
\begin{itemize}
\item for almost every $x\in \overline{{\rm X}(\mathbb{Z}[1/p])}
\subset {\rm X}(\mathbb{R})$, every $\delta>0$, and $\epsilon\in (0,\epsilon_0(x,\delta))$,
the system of inequalities
\begin{equation}\label{eq:ddioph}
\|x-r\|\le \epsilon\quad\hbox{and}\quad \hbox{\rm D}(r)\le 
\left(\epsilon^{-\frac{\dim({\rm X})}{a_p({\rm G})}-\delta}\right)^{q_p({\rm G})/2}
\end{equation}
has a solution $r\in {\rm X}(\mathbb{Z}[1/p])$.
\item for every $x\in {\rm X}(\mathbb{R})$, every $\delta>0$, and $\epsilon\in (0,\epsilon_0(\delta))$,
the system of inequalities
\begin{equation*}
\|x-r\|\le \epsilon\quad\hbox{and}\quad \hbox{\rm D}(r)\le 
\left(\epsilon^{-\frac{\dim({\rm X})}{a_p({\rm G})}-\delta}\right)^{q_p({\rm G})}
\end{equation*}
has a solution $r\in {\rm X}(\mathbb{Z}[1/p])$.
\end{itemize}
\end{thm}

In the paper \cite{ggn2} (joint with A.~Ghosh), we also prove analogues of Khnichin's and Jarnik's theorems
for Diophantine approximation on homogeneous varieties.

The quality of Diophantine approximation in Theorem \ref{th:dioph}
is controlled by the integrability exponent $q_p({\rm G})$. 
When $q_p({\rm G})=2$, then estimate (\ref{eq:ddioph}) gives
the best possible result up to arbitrary small $\delta>0$.
Coming back to the example (\ref{eq:sphere1}), we note that in this case
$\rm G$ is the orthogonal group and
$q_p({\rm G})=2$, which follows from the work of P.~Deligne \cite{de} on the classical
Ramanujan conjecture (see \cite{c} for a comprehensive treatment) combined with the
Jacquet--Langlands correspondence \cite{jl}.
This topic is surveyed in \cite{r}.

In conclusion we note the intriguing fact  that conversely, Diophantine approximation considerations
can provide insights in the direction of the generalised Ramanujan conjectures.
In particular, one can deduce estimates on the integrability exponents. For instance,
analysing the quality of Diophantine approximation in $\mathbb{R}^d$ by 
the rational points in $\mathbb{Z}[1/p]^d$, one can deduce that
$$
q_p(\hbox{SL}_d)\ge 2(d-1).
$$
We refer to \cite[Corollary~1.8]{ggn} for a more general result in this direction.

\section{Ergodic theorems for lattice subgroups}

\subsection{Ergodic theorems and equidistribution}
We turn to discuss ergodic theorems for countable groups such as $\Gamma=\hbox{SL}_n(\mathbb{Z})$
and more generally for arithmetic subgroups in semisimple algebraic groups.
When $G$ is a connected semisimple group, 
we have approached this problem via spectral methods, and detailed information about 
the unitary representation theory of $G$ was essential. However, for discrete lattice subgroups 
this approach is not feasible, since the irreducible unitary representations of a  lattice cannot
be classified (in a Borel manner), and a general unitary representation of a lattice cannot be represented uniquely as a direct integral of irreducible ones. Thus non-amenable discrete groups appear as the hardest case 
in ergodic theory because they possess neither an asymptotically invariant family of sets, nor a usable spectral theory, and so other approaches must be developed.  

Although our approach can be applied to lattices in general groups,
for the sake of keeping our exposition simple we will present it
for lattices $\Gamma$ in connected simple Lie groups $G$. For more general formulations, we refer to \cite{GN10}. 
We take a family of increasing compact domains $B_t$ in $G$, set $\Gamma_t=\Gamma\cap B_t$,
and denote by $\lambda_t$ the uniform probability measures supported on the sets $\Gamma_t$.
Given an arbitrary measure-preserving ergodic action of $\Gamma$ on a (standard) probability space $(X,\mu)$,
we consider the corresponding averaging operators
$$
\pi_X(\lambda_t) f(x)=\frac{1}{\abs{\Gamma_t}}\sum_{\gamma\in \Gamma_t} f(\gamma^{-1}x), \quad f\in L^p(X).
$$  
In order for the sets $\Gamma_t$ to behave in a regular manner, we assume that the domains $B_t$
are admissible in the sense of Definition \ref{def:addmiss}. This mild regularity assumption
turns out to be sufficient to establish general quantitative ergodic theorems:

\begin{thm}\label{lattices 2}
Assume that the action of  $\Gamma$ on $(X,\mu)$ has a spectral gap.  
Then the following holds. 
\begin{itemize}
\item \emph{Mean ergodic theorem :}
For every $f\in L^p(X)$ with  $ 1 < p < \infty$, 
$$
\norm{\pi_X(\lambda_t) f-\int_X f\, d\mu}_p \le C_p |\Gamma_t|^{-\theta_p}\norm{f}_p
$$
with uniform $C_p,\theta_p > 0$, independent of $f$. 
\item \emph{Pointwise ergodic theorem:}
For every $f\in L^p(X)$ with $p>1$ and  for $\mu$-almost every $x$, 
$$
\abs{\pi_X(\lambda_t)f-\int_X f\, d\mu}\le C_p(f,x)|\Gamma_t|^{-\theta'_p}
$$
with uniform $\theta_p'>0$, independent of $f$ and $x$.
\end{itemize}
\end{thm}

We emphasize that this result holds for all $\Gamma$-actions which have a spectral gap, and furthermore, if $\Gamma$ has property $T$ of Kazhdan, the parameters $C_p$, $\theta_p$ and $\theta_p^\prime$ are independent of the action. The only connection to the original
embedding of $\Gamma$ in the  group  $G$ is in the definition of the sets $\Gamma_t$. 
We remark that an $L^s$-norm bound for $C_p(f,\cdot)$ (for a suitable $s$ depending on $p$) holds. 
 Finally, when 
the action of $\Gamma$ does not have a spectral gap, we also establish mean and pointwise
ergodic theorem (but without a rate of convergence, of course). We refer to \cite{GN10} for details.
 
Theorem \ref{lattices 2} demonstrates a remarkable contrast with the ergodic theory of amenable groups
where no rate of convergence can be established, as discussed in Section 3. Let us give an example to illustrate
Theorem \ref{lattices 2}:
\begin{itemize}
\item Consider the action $\Gamma=\hbox{SL}_n(\ZZ)$ on the torus $\TT^n=\RR^n/\ZZ^n$
and set $$\Gamma_t=\{\gamma\in \Gamma:\, \log\|\gamma\|\le t\}.$$ Then
for every $f\in L^p(\TT^n)$ with $p > 1$ and almost every $x\in \TT^n$,
$$\abs{\frac{1}{\abs{\Gamma_t}}\sum_{\gamma\in \Gamma_t} 
f(\gamma^{-1}x)-\int_{\TT^n} f\, d\mu}\le C_p(f,x)|\Gamma_t|^{-\theta'_p }
$$
with uniform $\theta_p'>0$.
\end{itemize}
The orbit structure in the above example is quite complicated.
The torus contains the dense set of rational points which consists of 
finite $\Gamma$-orbits, but the orbit of every irrational point is 
dense. In particular, the size of the estimator $C_p(f,x)$ depends very 
sensitively on Diophantine properties of $x$, and it would be interesting
to make this dependence more explicit.  In the setting of random walks on the torus,
a similar problem was recently investigated in \cite{bflm}.

In the case of isometric actions of lattices,
we derive the following quantitative equidistribution result,
which again has no classical amenable analogue:

\begin{cor}\label{lattices 3}
Let $\Gamma$ act on isometrically on a compact Riemannian manifold $X$ equipped with a smooth measure
$\mu$ of full
support. If the action has spectral gap, than for every H\"older function $f\in C^a(X)$, and for every point $x$ in the manifold
$$
\frac{1}{|\Gamma_t|}\sum_{\gamma\in \Gamma_t} 
f(\gamma^{-1}x)=\int_{X} f\, d\mu+ O\left(|\Gamma_t|^{-\theta_a}\|f\|_{C^a}\right)
$$
with uniform $\theta_a > 0$. 
\end{cor}

\subsection{Method of the proof of the  ergodic theorems}

Our approach to proving the ergodic theorems for the averages $\pi_X(\lambda_t)$
is based on the idea of inducing the action of $\Gamma$  to obtain 
an action to $G$ on a larger space $Y=(G\times X)/\Gamma $. The latter space has the structure 
of a fiber bundle over the space $G/\Gamma$ with fibers isomorphic to $X$.
We then reduce the ergodic theorems for the averaging operators 
$\pi_X(\lambda_t)$, supported on $\Gamma\cap B_t$ and acting on $L^p(X)$, 
to the ergodic theorems for the averaging operators $\pi_Y(\beta_t)$,
supported on the domains $B_t$ and acting on $L^p(Y)$.
This involves developing a series of  approximation arguments which 
constitute generalisations of the argument used in the solution of the lattice point counting problem
in Section \ref{sec:counting}.
The essence of the matter is to approximate
$$
\pi_X(\lambda_t)\phi(x)=
\frac{1}{\abs{\Gamma\cap B_t}}\sum_{\gamma\in \Gamma\cap B_t} \phi(\gamma^{-1}x),\quad
\phi \in L^p(X),
$$
by
$$
\pi_Y(\beta_{t\pm c\vre})f_\vre(y)=\frac{1}{|B_{t\pm c\vre}|}\int_{B_{t\pm c\vre }} f_\vre(g^{-1}y)\, dg.
$$
with suitably chosen $f_\vre \in L^p(Y)$.
The link between the two expression is given by setting $y=(h,x)\Gamma$ and 
 $$f_\vre((h,x)\Gamma)=\sum_{\gamma\in \Gamma} \chi_\vre (h\gamma)\phi(\gamma^{-1}x),$$
where $\chi_\vre$ is the normalized characteristic function of an identity neighborhood ${O}_\vre$.

 The ergodic averages $\pi_Y(\beta_{t\pm c\vare})f_\vre$ can then be rewritten in full as 
$$\sum_{\gamma\in \Gamma}\left(\frac{1}{| B_{t\pm c\vare}|}\int_{B_{t\pm c\vare}}\chi_\vre(g^{-1}h\gamma)\,dg\right)\phi(\gamma^{-1}x),$$
We would like the expression in parentheses to be equal to $1$ when $\gamma\in \Gamma\cap B_{t-c\vare}$ and equal to $0$ when  $\gamma\notin \Gamma\cap B_{t+c\vare}$, in order to be able to compare it to $\pi_X(\lambda_t)\phi$.  
These favourable lower and upper estimates depend only on the regularity properties of the sets $B_t$, namely,
on the admissibility property.
Indeed, the lower bound arises since if  $\chi_\vre(g^{-1}h\gamma)\neq 0$  and $g\in B_{t-c\vare}$,
then $\gamma\in B_{t}$. The upper bound holds since for $\gamma\in \Gamma\cap B_t$,  the
support of function $g\mapsto \chi_\vre(g^{-1}h\gamma)$ is contained in $B_{t+c\vare }$. 

A fundamental point in the completion of the proof of the ergodic theorems is played by
an \emph{invariance principle}, namely by the fact that for any given $f\in L^p(Y)$, 
the pointwise ergodic theorem for the averages $\pi_Y(\beta_t)$ holds for a set of points which contain a
$G$-invariant set of full measure.
 Since the space $Y=(G\times X)/\Gamma$ is a $G$-equivariant bundle over the $G$-transitive space $G/\Gamma$,  
this implies that for \emph{every point} $h\Gamma$, the set of points $x\in X$ in the fiber above it where 
the pointwise ergodic theorem holds is conull in $X$. This allows us to deduce that the set of points
where $\pi_X(\lambda_t)\phi(x)$ converges is also a set of full measure in $X$.

\section{Distribution of orbits on algebraic varieties}
\subsection{Distribution of orbits in the de-Sitter space} 
Let us now turn to discuss the distribution of dense orbits for actions of arithmetic groups on algebraic varieties.
These action do not  preserve a finite measure (besides the measures supported on periodic points),
so that this problem falls in the realm of infinite ergodic theory.

Consider an action of a discrete group $\Gamma$ on a space $X$.
The distribution of orbits of $\Gamma$ is
encoded by the properties of the averaging operators
\begin{equation}\label{eq:average}
\mathcal{A}_t\phi(x)=\frac{1}{|\Gamma_t|}\sum_{\gamma\in\Gamma_t} \phi(\gamma^{-1}x)
\end{equation}
defined for an increasing sequence of finite subsets of $\Gamma$ that exhaust $\Gamma$
and for a family of test functions $\phi$ on $X$.
For example, if $\Gamma=\left<\gamma\right>$ is a cyclic group,  
and the action preserves a (possibly infinite) measure $\mu$, then
as in Section 1 we obtain that
\begin{equation}\label{eq:erg} 
\lim_{N\to\infty}\frac{1}{2N+1}\sum_{n=-N}^N \phi(\gamma^n x)\to \mathcal{P}\phi(x)\quad\hbox{in $L^2(X)$,}
\end{equation}
 where $\mathcal{P}$ denotes the orthogonal projection from
$L^2(X)$ to the subspace of $\Gamma$-invariant functions in $L^2(X)$. However, if the measure is infinite
and the action is ergodic, then
$\mathcal{P}=0$, and (\ref{eq:erg}) contains no information about the distribution of orbits.
Initially, one might expect that if we pick a correct normalisation factor such that $\frac{a(N)}{N}\to 0$,
then the averages $\frac{1}{a(N)}\sum_{n=-N}^N \phi(\gamma^n x)$ exhibit a nontrivial 
limit. But in fact, J. Aaranson \cite[\S2.4]{aa} showed that if the action of 
the cyclic group $\Gamma=\left<\gamma\right>$ is ergodic and conservative, then for
any normalisation $a(N)$, either for every nonnegative $\phi\in L^1(X)$,
$$
\liminf_{N\to\infty} \frac{1}{a(N)}\sum_{n=-N}^N \phi(\gamma^n x)=0\quad\hbox{almost everywhere,}
$$
or there exists a subsequence $N_k$ such that
such that for every nonnegative $\phi\in L^1(X)$, $\phi\ne 0$,
$$
\lim_{k\to\infty} \frac{1}{a(N_k)}\sum_{n=-N_k}^{N_k} \phi(\gamma^n x)=\infty\quad\hbox{almost everywhere.}
$$

But while the averages for orbits of $\mathbb{Z}$-actions on infinite measure spaces
tend fluctuate wildly, it turns out that orbits of ``large'' groups behave in
a more regular fashion. 
We will presently demonstrate this point by computing the distribution of the orbits of the groups of isometries
of the $d$-dimensional de Sitter space $dS_d$. This space can be realised
as a hypersurface 
$$
x_1^2+\cdots + x_d^2-x_{d+1}^2=1
$$
in $\mathbb{R}^{d+1}$ equipped with the Minkowski metric $ds^2=dx_1^2+\cdots+dx_d^2-dx_{d+1}^2$.
The problem of distribution of orbits on the de Sitter space was raised by V. Arnol'd 
\cite[1996-15, 2002-16]{a}.
After stating the result,  we will  explain a 
general strategy how to study distribution of orbits on homogeneous spaces.
Our exposition is based on \cite{gw,gn_duality}. We also refer to 
\cite{le,n1,m,lp,g1,g2,gm,gui,n2,mw}
for related results.

\begin{figure}[h]
\includegraphics[height=1.9in]{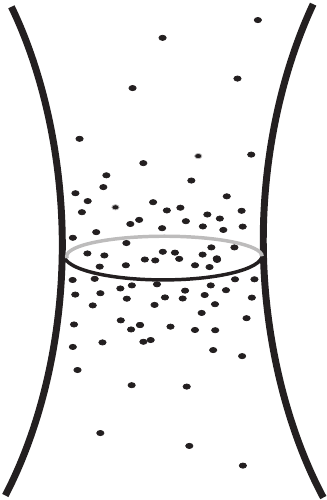} 
\caption{The de Sitter space}
\end{figure}

In the setting of Arnol'd's question
regarding the de Sitter space, we prove the following.
Let $\Gamma$ be a discrete subgroup of $\hbox{Isom}(dS_d)$ with finite covolume.
We identify $\hbox{Isom}(dS_d)$ with the orthogonal group $\hbox{SO}(d,1)$ and
denote by $\|\cdot\|$ the standard Euclidean norm on the space of $(d+1)$-dimensional
matrices.

\begin{itemize}
\item When $d=2$, for every $\phi\in L^1(dS^d)$ with compact support
\begin{equation}\label{eq:ds_eq1}
\lim_{t\to\infty} \frac{1}{t}\sum_{\gamma\in\Gamma:\, \log \|\gamma\|\le t}
 \phi(\gamma^{-1}x)= \int_{dS^d} \phi\, d\nu\quad\hbox{almost everywhere},
\end{equation}
where $\nu$ denotes a (nonzero) invariant measure on $dS_d$ (whose normalization depends on $\Gamma$).
\item When $d\ge 3$, for every continuous function $\phi$ with compact support
and every $x\in dS^d$ with dense $\Gamma$-orbit,
\begin{equation}\label{eq:ds_eq2}
\lim_{t\to\infty}\frac{1}{e^{(d-2)t}}\sum_{\gamma\in\Gamma:\, \log \|\gamma\|\le t}
 \phi(\gamma^{-1}x)= \int_{dS^d} \phi\, d\nu_x,
\end{equation}
where the limit measure is $d\nu_x(y)=\frac{d\nu(y)}{(1+\|x\|^2)^{(d-2)/2}(1+\|y\|^2)^{(d-2)/2}}$.

Moreover, if in addition $\phi$ is subanalytic and Sobolev, then for almost every $x\in dS^d$,
\begin{equation}\label{eq:ds_eq3}
\frac{1}{e^{(d-2)t}}\sum_{\gamma\in\Gamma:\, \log \|\gamma\|\le t}
 \phi(\gamma^{-1}x)= \int_{dS^d} \phi(y)\, d\nu_x(y)+O_{\phi,x}(e^{-\delta t})
\end{equation}
with $\delta=\delta(\phi,x)>0$.
\end{itemize}

The above equidistribution results have several remarkable features. In the case of dimension $2$, 
while the cardinality of the set $\{\gamma\in\Gamma:\, \log \|\gamma\|\le t\}$ 
grows exponentially as $\hbox{const}\cdot e^t$, it turns out that only for  a polynomial number of points $\gamma$ does the orbit points $\gamma x$ 
in (\ref{eq:ds_eq1}) come back to a compact set. Nonetheless, the small fraction of points returning by time $t$ 
 becomes equidistributed in the compact set, with respect to the invariant measure as $t\to \infty$.
 
When $d\ge 3$, while the cardinality of the set $\{\gamma\in\Gamma:\, \log \|\gamma\|\le t\}$ grows as  $\hbox{const}\cdot e^{(d-1)t}$, only an exponentially small fraction of them satisfy that $\gamma x$ returns to a compact set. The set of returning points does become  equidistributied in the compact set, but this time the limiting measures $\nu_x$
are not invariant under $\Gamma$, and furthermore depend nontrivially on $x$.
The measures $\nu_x$, $x\in X$, on the pseudo-Riemannian space $X$
should be considered as analogues of the Patterson--Sullivan measures in this context.

\subsection{The duality principle on homogeneous spaces}
The techniques developed in \cite{gn_duality} apply in the following general setting.
Let $G\subset \hbox{SL}_d(\mathbb{R})$ be  connected and of finite index in an  algebraic group,  and 
$\Gamma$ a discrete subgroup of $G$ with finite covolume.
Let $X$ be an algebraic homogeneous space of $G$ equipped with a smooth measure
on which $\Gamma$ acts ergodically.
We fix a proper homogeneous polynomial $p:\hbox{M}_d(\mathbb{R})\to [0,\infty)$
and consider the sets 
$$\Gamma_t=\{\gamma\in\Gamma:\, \log p(\gamma)\le t\}.$$
Our goal is to describe the asymptotic distribution of orbits of $\Gamma$,
or in other words the asymptotic behaviour of the sums
$\sum_{\gamma\in\Gamma_t} \phi(\gamma^{-1}x)$ as $t\to\infty$
for a sufficiently rich collection of functions $\phi$ on $X$ with compact support.
In order to find the right normalisation for the sum
$\sum_{\gamma\in\Gamma_t} \phi(\gamma^{-1}x)$,
one needs to compute what proportion of points in the orbit
returns to compact subsets $\Omega\subset X$. Given a compact subset $\Omega$ of $X$
with non-empty interior and $x\in X$, we set
$$
a=\limsup_{t\to \infty} \frac{\log |\Gamma_t^{-1} x\cap \Omega|}{t}.
$$
One can check that this quantity is independent of the choices of $\Omega$ and $x$.
We will distinguish two cases according to whether $a=0$, i.e., the return rates are
at most subexponentially, and $a>0$, i.e., the return rates are exponential.

\begin{thm}[polynomial return rates]\label{th:polynomial}
Assume that $a=0$.
Then there exists $b\in \mathbb{Z}_{>0}$ such that the averaging operator
$$
\mathcal{A}_t\phi(x)=\frac{1}{t^b}\sum_{\gamma\in \Gamma_t}
\phi(\gamma^{-1}  x)
$$
satisfy
\begin{itemize}
\item {\rm (strong maximal inequality)} For every $p>1$, a compact $D\subset X$,
and $\phi\in L^p(D)$,
$$
\left\|\sup_{t\ge t_0} |\mathcal{A}_t\phi|\right\|_{L^p(D)}\le \hbox{\rm const}(p,D)\cdot \|\phi\|_{L^p(D)}
$$

\item {\rm (pointwise ergodic theorem)}
For every $\phi\in L^1(X)$ with compact support,
$$
\lim_{t\to\infty} \mathcal{A}_t\phi(x)=\int_X\phi\,d\nu\quad\hbox{almost everywhere},
$$
where $\nu$ is a (nonzero) $G$-invariant measure on $X$.
\end{itemize}
\end{thm}

When $a>0$, we prove analogous results for functions lying in a suitable Sobolev space $L^p_l(X)$.
In this case it is also possible to establish rates of convergences.

\begin{thm}[exponential return rates]\label{th:exponential}
Assume that $a>0$.
Then there exist  $b\in \mathbb{Z}_{\ge 0}$ and $l\in \mathbb{Z}_{\ge 0}$
such that the averaging operator
$$
\mathcal{A}_t\phi(x)=\frac{1}{e^{at}t^b}\sum_{\gamma\in \Gamma_t}
\phi(\gamma^{-1}\cdot  x)
$$
satisfies
\begin{itemize}
\item {\rm (strong maximal inequality)} For every $p>1$, a compact $D\subset X$,
and a nonnegative $\phi\in L^p_l(D)$,
$$
\left\|\sum_{t\ge t_0} |\mathcal{A}_t\phi|\right\|_{L^p(D)}\le \hbox{\rm const}(p,D)\cdot \|\phi\|_{L_l^p(D)}
$$

\item {\rm (pointwise ergodic theorem)}
For every $p>1$ and a nonnegative bounded $\phi\in L^p_l(X)$ with compact support,
\begin{equation}\label{eq:exp2}
\lim_{t\to\infty} \mathcal{A}_t\phi(x)=\int_X\phi\,d\nu_x\quad\hbox{almost everywhere},
\end{equation}
where $\nu_x$, $x\in X$, is a family of absolutely continuous measures on $X$ with
positive densities.

\item {\rm (quantitative pointwise ergodic theorem)}
For every $p>1$ and  a nonnegative continuous subanalytic $\phi\in L^p_l(X)$ with compact support,
\begin{equation}\label{eq:exp3}
\mathcal{A}_t\phi(x)=\int_X\phi\,d\nu_x
+\sum_{i=1}^b c_i(x,\phi)t^{-i} +O_{\phi,x} \left(e^{-\delta t}\right)
\quad\hbox{almost everywhere}
\end{equation}
with $\delta=\delta(x,\phi)>0$.
\end{itemize}
\end{thm}

If in Theorem \ref{th:exponential} we additionally assume that the stabiliser of  a point $x$ in $G$
is semisimple, then the stated results can be further improved. One can replace the Sobolev norms
by $L^p$ norms (see \cite[Theorem~1.3]{gn_duality}). Moreover, 
one can show (see \cite{gw}) for continuous functions $\phi$ with compact support, 
the asymptotic formula (\ref{eq:exp2}) holds for all $x\in X$ whose $\Gamma$-orbit in $X$ is dense.

\subsection{Equidistribution in algebraic number fields}

It seems unlikely that in the generality of Theorem \ref{th:exponential}, the asymptotic expansion
(\ref{eq:exp3}) can be proved for an  explicit set of $x\in X$ of positive measure,
but there is an important case when this can be achieved. 
Assuming that the action of $\Gamma$ on $X$ preserves a Riemannian metric, one can show
(see \cite{gn_isometry}) that all $\Gamma$-orbits in $X$ become equidistributed with 
a rate. The observation that all orbits for isometric actions on compact spaces
are equidistributed goes back to the
work of Y.~Guivarc'h \cite{gu}, and in \cite{gn_isometry} we have developed a quantitative
version of this argument which also applies to infinite-measure setting.
For instance, fix a non-square integer $d >0$ and let us consider the action of $\Gamma=\hbox{SL}_2(\mathbb{Z}[\sqrt{d}])$ on 
the upper half plane $\mathbb{H}^2$ by fractional linear transformations.
We set 
$$
\Gamma_t=\{\gamma\in\Gamma:\, \log(\|\gamma\|^2+\|\bar \gamma\|^2)\le t\}
$$
where $\|\cdot \|$ is the standard Euclidean norm, 
and $\bar\gamma$ denotes the Galois involution of the field $\mathbb{Q}(\sqrt{d})$. 
Since $\Gamma$ embeds diagonally in $G=\hbox{SL}_2(\mathbb{R})\times \hbox{SL}_2(\mathbb{R})$
as an irreducible lattice, and $\mathbb{H}^2$ is a homogeneous space of the group $G$,
Theorem \ref{th:exponential} applies in this case and gives information about
the asymptotic distribution of orbits $\Gamma x$ for almost every $x\in X$, and
the method of \cite[Theorem~1.5]{gn_isometry} leads to the following equidistribution result:
there exists $\delta>0$ such that given a compact subset $D\subset \mathbb{H}^2$, 
for every $x\in  D$ and every $\phi\in C^1(\mathbb{H}^2)$ with support contained in $D$,
$$
\frac{1}{e^t}\sum_{\gamma\in\Gamma_t} \phi(x\gamma)=\int_{\mathbb{H}^2} \phi(z)\, d\nu(z)
+O_D(e^{-\delta t}),
$$
where $\nu$ is a (nonzero) $G$-invariant measure on $\mathbb{H}^2$.

\subsection{Ideas of the proof}
In conclusion we indicate some ideas that lie behind the proof of Theorem \ref{th:polynomial}
and \ref{th:exponential}. The argument can be divided into two main steps:
\begin{enumerate}
\item[(I)] compare the asymptotic behaviour of discrete averages 
\begin{equation}\label{eq:discrete}
\sum_{\gamma\in \Gamma:\,\log p(\gamma)\le t} \phi(\gamma^{-1}x)
\end{equation}
with the asymptotic behavior of continuous averages
\begin{equation}\label{eq:continuous}
\int_{g\in G:\,\log p(g)\le t} \phi(g^{-1}x)\, dm_G(g),
\end{equation}
where $m_G$ is an invariant measure on $G$.
  
\item[(II)] establish the asymptotics for the continuous averages (\ref{eq:continuous}).
\end{enumerate}

To achieve step (I), we think of $G$ as a fiber bundle over $X$, and show 
that the fibers of this bundle become equidistributed 
on the space $\Gamma\backslash G$. More explicitly, we identify
$X$ with the homogeneous space $G/H$ where $H$ is a subgroup of $G$. For $u,v\in G$, we set
$$
H_t[u,v]=\{h\in H:\, \log p(uhv)\le t\}.
$$
Then we show that the sum (\ref{eq:discrete}) can be approximated 
by a finite linear combination of integrals
$$
\int_{H_{t_i}[u_i,v_i]} f_i(y_ih)\,dm_H(h)
$$
with suitably chosen $t_i\approx t$, $u_i,v_i\in G$, $f_i:\Gamma\backslash G\to \mathbb{R}$, and $y_i\in
\Gamma\backslash G$,
where $m_H$ is a left invariant measure on $H$.
The asymptotic behaviour of these integrals can be analysed using 
either representation-theoretic techniques as in \cite{gn_duality} or
the theory of unipotent flows as in \cite{gw}. In both cases we conclude that
as $t\to\infty$,
$$
\frac{1}{|H_t[u,v]|} \int_{H_t [u,v]} f(y h)dm_H(h) 
\to \int_{\Gamma\backslash G} f\, dm_{\Gamma\backslash G},
$$
where $m_{\Gamma\backslash G}$ denotes the normalised invariant measure on ${\Gamma\backslash G}$.
Finally, this allows to conclude that the above linear combination is approximately equal to
the integral (\ref{eq:continuous}), which completes step (I).

In step (II), we establish an asymptotic formula for integral
$$
v(t)=\int_{g\in G:\,\log p(g)\le t} \phi(g^{-1}x)\, dm(g)
$$
where $\phi$ is a continuous subanalytic function with compact support.
For this, we consider its transform
$$
f(s)=\int_0^\infty t^{-s} v(\log t)\,dt
$$
which converges when $\hbox{Re}(s)$ is sufficiently large.
Using a suitable version of resolution of singularities as in \cite{p}
we show that $f(s)$ has a meromorphic continuation beyond the first pole.
Then the asymptotic formula for $v(t)$ follows via a Tauberian 
argument.

\section{Acknowledgements}
Part of this paper was written in Spring 2011 when A.G. was visiting 
\'Ecole Polytechnique F\'ed\'erale de Lausanne
during the programme ``Group Actions in Number Theory'', and 
he would like to thank Emmanuel Kowalski, Philippe Michel,
and the Centre Interfacultaire Bernoulli 
for hospitality.

Both authors would like to express their gratitude
to the Ergodic Theory Group at the F\'ed\'eration Denis Poisson for  the 
opportunity to explain the present work in the lecture series ``Th\'eorie ergodique
des actions de groupes'' held at the University of Tours in April, 2011. In particular we would like to thank the organizers, Claire Anantharaman-Delaroche, Jean-Philippe Anker, 
and Emmanuel Lesigne. 

A.G. also would like to thank for hospitality Manfred Einsiedler 
and the Institute for Mathematical Research at ETH, Z\"urich 
during his visit in Autumn 2012 when the work on this survey was completed.

A.G. was supported in part by EPSRC, ERC, and RCUK, and 
A.N. was supported by an ISF grant.

\end{document}